\documentclass[reqno]{amsart}

 \usepackage{mathtools}
 \usepackage{amsthm,amssymb,mathrsfs,enumitem,wasysym,hyperref}
 \usepackage{esint,enumitem}

 \usepackage[numbers,sort&compress]{natbib}
 \usepackage{mathscinet}

 \numberwithin{equation}{section}
 \numberwithin{figure}{section}
 \theoremstyle{plain}
 \newtheorem{thm}{Theorem}[section]
   \theoremstyle{definition}
   
   \newtheorem*{defn*}{Definition}
   \theoremstyle{remark}
   \newtheorem*{rem*}{Remark}
   \theoremstyle{plain}
   
   \theoremstyle{plain}
   \newtheorem{prop}[thm]{Proposition}
   
   \theoremstyle{definition}

   \newcounter{step}

 \makeatletter
 \newcommand{\norm}{\@ifstar{\@normb}{\@normi}}
 \newcommand{\@normb}[2]{\left\Vert{#1}\right\Vert_{#2}}
 \newcommand{\@normi}[2]{\Vert{#1}\Vert_{#2}}
 \makeatother

 \global\long\def\tSob#1{{H}^{#1}} 
 \global\long\def\thSob#1{\dot{H}^{#1}}

 \global\long\def\Leb#1{L^{#1}} 
 \global\long\def\Lebdiv#1{L^{#1}_\sigma}

 \RequirePackage{bm}
 
 \newcommand{\boldb}{\bm{b}}
 \newcommand{\boldu}{\bm{u}}

 \DeclareMathOperator{\Div}{div}
 \newcommand{\relphantom}[1]{\mathrel{\phantom{#1}}}
 
 \newcommand{\myd}[1]{\,d{#1}}

\newcommand{\pbar}{\overline{p}}

\makeatletter
\def\@tocline#1#2#3#4#5#6#7{\relax
  \ifnum #1>\c@tocdepth % then omit
  \else
    \par \addpenalty\@secpenalty\addvspace{#2}%
    \begingroup \hyphenpenalty\@M
    \@ifempty{#4}{%
      \@tempdima\csname r@tocindent\number#1\endcsname\relax
    }{%
      \@tempdima#4\relax
    }%
    \parindent\z@ \leftskip#3\relax \advance\leftskip\@tempdima\relax
    \rightskip\@pnumwidth plus4em \parfillskip-\@pnumwidth
    #5\leavevmode\hskip-\@tempdima
      \ifcase #1
       \or\or \hskip 1em \or \hskip 2em \else \hskip 3em \fi%
      #6\nobreak\relax
    \hfill\hbox to\@pnumwidth{\@tocpagenum{#7}}\par% <---- \dotfill -> \hfill
    \nobreak
    \endgroup
  \fi}
\makeatother

\newcommand{\Proj}{\mathbb{P}}

\newcommand{\boldJ}{\mathbf{J}}
\newcommand{\boldE}{\mathbf{E}}

\DeclareMathOperator{\dv}{div}

\newcommand{\eqn}{\begin{eqnarray}}
\newcommand{\een}{\end{eqnarray}}

\usepackage[margin=1.2in]{geometry}

\usepackage{changes}

\definechangesauthor[name=Kwon, color=blue]{K}
\definechangesauthor[name=Shin, color=blue]{S}

 \title[Stokes-Magneto equations]{Global solutions to Stokes-Magneto equations with fractional dissipations}
 
\author{Hantaek Bae}
\address{Department of Mathematical Sciences, Ulsan National Institute of Science and Technology (UNIST), Republic of Korea}
\email{hantaek@unist.ac.kr} 

\author{Hyunwoo Kwon}
\address{Division of Applied Mathematics, Brown University, Providence, Rhode Island, USA}
\email{hyunwoo\_kwon@brown.edu} 

\author{Jaeyong Shin}
\address{Department of Mathematical Sciences, Ulsan National Institute of Science and Technology (UNIST), Republic of Korea}
\email{sinjaey@unist.ac.kr}

\subjclass{35Q35;76W05}
\keywords{Global existence; Uniqueness; Strong solutions; Mild solutions; Fractional diffusions}

 \allowdisplaybreaks

\begin{document}

\begin{abstract}
In this paper, we investigate a Stokes-Magneto system with fractional diffusions. We first deal with the non-resistive case in $\mathbb{T}^{d}$ and establish the local and global well-posedness with initial magnetic field $\boldb_0\in H^{s}(\mathbb{T}^d)$. We also show the existence of a unique mild solution of the resistive case with initial data $\boldb_0$ in the critical $L^{p}(\mathbb{R}^d)$ space. Moreover, we show that $\|\boldb(t)\|_{L^{p}}$ converges to zero as $t\rightarrow\infty$ when the initial data is sufficiently small.
\end{abstract}

\maketitle

\vspace{-3ex}

%%%%%%%%%%%%
\section{Introduction}
 %%%%%%%%%%%

%%%%%%%%%%%%%%%%%%%%%%%%%%%%%
\subsection{Motivation: Magnetic relaxation}
%%%%%%%%%%%%%%%%%%%%%%%%%%%%
The most common description of the electromagnetic field uses two vector fields: the electric field $\boldE$ and the magnetic field $\boldb$. The electromagnetic field is governed by Maxwell's equations:
\begin{subequations} \label{MHD}
\begin{align}
\text{Ampere's Law:} \quad & \nabla\times \boldb=\boldJ, \label{MHD c}\\
\text{Faraday's Law:} \quad &  \nabla\times \boldE=-\partial_t \boldb,\label{MHD d}\\
\text{Ohm's Law:} \quad &\boldE+\boldu\times \boldb=\eta \boldJ, \label{MHD e}\\
\text{Incompressibility:} \quad& \dv \boldb=0, \label{MHD f}
\end{align}
\end{subequations}
where $\boldu$ is the velocity field of electrically conducting fluid, $\boldJ$ is the current density, and $\eta$ is the resistivity (or magnetic diffusivity) of $\boldb$. From the above equations, we derive 
\begin{equation}\label{Eq of B}
\partial_t \boldb- \eta \Delta \boldb+(\boldu\cdot \nabla) \boldb -(\boldb\cdot \nabla) \boldu =0.
\end{equation}
 This is the so-called induction equation of the magnetohydrodynamics of $\boldb$ if $\boldu$ is known (see e.g.  \cite{D01} for a detailed derivation of the equation).
 
 Due to the mathematical interest and physical importance, such as nuclear fusions (see e.g. \cite{Arn74,G67} and references therein), much attention has been paid to the way of constructing magneto-hydrostatic equilibrium. Related to our results, we begin with the case $\eta=0$ so that the magnetic helicity is conserved under (\ref{Eq of B}). We also take the spatial domain to be $\mathbb{T}^{d}$. When a conducting fluid, initially at rest, is permeated by a magnetic field, the associated Lorentz force $\mathbb{F}=\boldJ\times \boldb$ drives motion in the fluid:
\begin{equation}  \label{Eq of u}
-\nu\Delta \boldu+\nabla p_0=\boldJ\times \boldb,
\end{equation}
where $\nu$ is the viscosity of the fluid and $p_0$ is the fluid pressure. If we write the magnetic energy by $M(t)=\frac{1}{2} \|\boldb(t)\|^{2}_{L^{2}(\mathbb{T}^{d})}$ and the magnetic helicity by $\mathcal{H}(t)$, then the Arnold inequality \cite{Arn74} on $\mathbb{T}^{d}$  gives a positive lower bound of $M(t)$:
\[
M(t)\ge {|\mathcal{H}(t)|=|\mathcal{H}(0)|}>0
\]
when the initial helicity is {non-zero}.  Moreover, from (\ref{Eq of B}) and (\ref{Eq of u}), we have
\[
\frac{d}{dt}M(t) =-\nu\left\|\nabla \boldu\right\|^{2}_{L^{2}(\mathbb{T}^d)}.
\]
Provided that the smooth solutions of (\ref{Eq of B}) and (\ref{Eq of u}) exist globally-in-time, the minimum energy is attained only when the fluid comes to rest again. In this situation, 
\[
\boldu(t,x)\rightarrow 0, \quad \boldb(t,x)\rightarrow \overline{\boldb}(x)
\]
as $t\rightarrow \infty$ and for all $x$.  This process is called the magnetic relaxation introduced by Moffatt \cite{M85, M21}. Therefore, it is necessary to find a global-in-time solution of (\ref{Eq of B}) and (\ref{Eq of u}) as a precondition for verifying the magnetic relaxation. We also remark that Constantin-Pasqualotto \cite{CP22} constructed {magneto-hydrostatic} equilibria as long time limits of certain Voigt regularizations of the MHD equations. See also \cite{ELP21,BL96} for different constructions of magneto-hydrostatic equilibrium.

%%%%%%%%%%%%%%%%%%%%%%%%%%%%%%%%%
\subsection{Stokes-Magneto system with fractional dissipations}
%%%%%%%%%%%%%%%%%%%%%%%%%%%%%%%%%%
Toward a rigorous theory to reach magnetic relaxation, we investigate the following Stokes-Magneto system  with fractional dissipations:
\begin{equation}\label{eq:MRE}
\left\{
\begin{alignedat}{2}
\nu \Lambda^{2\alpha} \boldu +\nabla \pbar &=(\boldb\cdot \nabla)\boldb&&\quad \text{in } \Omega\times (0,T),\\
\partial_t \boldb+\eta \Lambda^{2\beta} \boldb +(\boldu\cdot \nabla)\boldb &=(\boldb\cdot \nabla)\boldu &&\quad \text{in }\Omega \times(0,T),\\
\Div \boldu =\Div \boldb &=0 &&\quad \text{in } \Omega\times (0,T),\\
\boldb(\cdot,0)&=\boldb_0 &&\quad \text{on } \Omega.
\end{alignedat}\right.
\end{equation}
Here $\Omega$ is either $\mathbb{R}^d$ or $\mathbb{T}^d$, $\boldu=(u^1,u^2,\dots,u^d):\Omega\times [0,T)\rightarrow \mathbb{R}^d$ is the fluid velocity field, $\boldb=(b^1,b^2,\dots,b^d): \Omega\times [0,T)\rightarrow \mathbb{R}^d$ is the magnetic field, and $\pbar:\Omega\times [0,T)\rightarrow \mathbb{R}$ denotes the total pressure $\pbar = p_0 + \frac{1}{2}|\boldb|^2$, with $p_0$ being the pressure on the fluid. The nonnegative constants $\nu$ and $\eta$ stand for the viscosity constant and the magnetic diffusivity, respectively. We only consider the case $\nu>0$, so set $\nu=1$ for simplicity. Here $\Lambda^s$ are fractional Laplacians defined via the Fourier transform in (\ref{FL1}) on $\mathbb{R}^d$ and (\ref{FL2}) on $\mathbb{T}^d$, respectively.  We note that (\ref{eq:MRE}) is exactly the model introduced by Moffatt \cite{M21} when $\alpha=1$ and $\eta=0$. 

\vspace{1ex}

Unlike classical magnetohydrodynamics equations (MHD in short) (see \cite{Feng23} and references therein for a detailed list of recent results for MHD), there are only a few mathematical results  for \eqref{eq:MRE} and we list some of them separately for the case $\eta=0$ and $\eta>0$.

\subsubsection*{The case $\eta=0$} When $\alpha=1$, the local existence and uniqueness of strong solutions of \eqref{eq:MRE} in $\Omega=\mathbb{R}^d$, $d=2,3$,  were established by Fefferman-McCormick-Robinson-Rodrigo \cite{FMRR14}. Beekie-Friedlander-Vicol \cite{BFV22} proved the local well-posedness of strong solutions to \eqref{eq:MRE}  in $\Omega=\mathbb{T}^d$, $d=2,3$, when $\alpha \geq 0$. They also proved that the local solution becomes global when $\alpha>\frac{d}{2}+1$, $d=2,3$. Meanwhile, Brenier \cite{Bre14} proved the global existence of dissipative weak solutions in $\Omega=\mathbb{T}^2$ when $\alpha=0$.
\subsubsection*{The case $\eta>0$} When $\alpha=\beta=1$, McCormick-Robinson-Rodrigo \cite{MRR14} proved the existence of weak solutions in various domains in 2 and 3 dimensions (uniqueness in 2 dimensions also). Later, the long-time behavior of weak solutions was investigated by Ji-Tan \cite{JT21-1} in $\Omega=\mathbb{R}^3$. When $\alpha=\beta=1$ and $\Omega=\mathbb{R}^3$, several blow-up conditions were derived by Tan \cite{Tan23} (which also contains a global-in-time solution with the smallness condition in Besov spaces). When $\beta\ge \frac{3}{2}$, Ji-Tan \cite{JT21} proved the global existence of strong solutions in $\Omega=\mathbb{R}^3$. Very recently, the global existence of weak solution in $\Omega=\mathbb{R}^d$ was established  by Kim and the second author \cite{KK23} when 
\begin{align}
\frac{1}{2}<\alpha<\frac{d+1}{2}, \quad &\beta>0, \quad \min\{\alpha+\beta, 2\alpha+\beta-1\}>\frac{d}{2}\nonumber
\intertext{which can be unique if}
\frac{1}{2}<\alpha<\frac{d+1}{2}, \quad &\beta \ge 1, \quad \min\{\alpha+\beta, 2\alpha+\beta-1\}\geq \frac{d}{2}+1.\label{eq:uniqueness-weak}
\end{align}

The goal of this paper is twofold.
 \begin{enumerate}
\item Improvement of  \cite{BFV22}: (i) weakening of the range of $s$ of the initial data space $\boldb_0 \in H^{s}$ and the upper bound of solutions, (ii)  providing a sufficient condition for the continuation in time of solution. In particular, we prove that \eqref{eq:MRE} is globally well-posed even when $\alpha=\frac{d}{2}+1$.
\item Existence of mild solution to \eqref{eq:MRE} in $\mathbb{R}^d$ with $\boldb_{0} \in L^{p}$ and  $(\alpha,\beta)$ satisfying
\[
\frac{1}{2}<\alpha< \frac{d+1}{2},\quad \beta>\frac{1}{2},\quad \alpha+\beta < d+1,\quad p=\frac{d}{\alpha+\beta-1}
\]
and 
\[
3\alpha+\beta<\frac{3d+4}{2}
\]
in addition if 
\[
\frac{d+2}{4}<\alpha<\frac{d+1}{2}.
\]
\end{enumerate}

{The remaining of this paper proceeds as follows. In Section \ref{sec:2}, we introduce our notations and the main results of this paper. Section \ref{sec:prelim} is devoted to giving and proving some preliminary results including logarithmic Sobolev inequality and estimates on fractional heat operators. The proof of main theorems will be given in  Sections \ref{sec:4} and \ref{sec:5}.}

%%%%%%%%%%%%%%%%%%%%%%%%%%%
\section{Notations and Main results}\label{sec:2}
%%%%%%%%%%%%%%%%%%%%%%%%%%%

%%%%%%%%%%%%%%%%%%%%%
\subsection{Notation and definitions}
%%%%%%%%%%%%%%%%%%%%%
All  constants are denoted by $C$ and we follow the convention that such constants can vary from expression to expression and even between two occurrences within the same expression. We write $C(p_{1}, p_{2}, \cdots, p_{n})$ to mean a constant that depends only upon the parameters $p_{1}, p_{2}, \cdots, p_{n}$. We say $A\apprle B$ if there exists a constant $C>0$ such that $A \leq CB$, where $C$ does not depend on $A$ and $B$. We also write $A\approx B$ if $A\apprle B$ and $B\apprle A$.

\vspace{1ex}

Let $\mathbb{R}^d$ denote $d$-dimensional Euclidean spaces of the points and let $\mathbb{T}^d=[0,2\pi]^d$ be $d$-dimensional torus. Let $C^\infty(\mathbb{R}^d)$ be the space of all smooth functions and $C^\infty(\mathbb{T}^d)$ be the space of all functions $C^\infty(\mathbb{R}^d)$ that are $2\pi$-periodic. We write $C^\infty_c(\mathbb{R}^d)$ and $\mathcal{S}(\mathbb{R}^d)$  the space of all smooth functions $C^\infty(\mathbb{R}^d)$ which have compact support in $\mathbb{R}^d$ and Schwartz class on $\mathbb{R}^d$, respectively. 

For $1\leq p<\infty$ and a function $f:\mathbb{R}^d\rightarrow \mathbb{R}$, let 
\[
\norm{f}{\Leb{p}(\mathbb{R}^d)}=\left(\int_{\mathbb{R}^d} |f|^p \myd{x}\right)^{1/p}
\]
be the $\Leb{p}$-norm of $f$. For $2\pi$-periodic function $f:\mathbb{T}^d\rightarrow \mathbb{R}$, we can also define $\norm{f}{\Leb{p}(\mathbb{T}^d)}$. Similarly, we can define the $\Leb{\infty}$-norm on $\mathbb{R}^d$ and $\mathbb{T}^d$. Note that $C_c^\infty(\mathbb{T}^d)=C^\infty(\mathbb{T}^d)$. {For $\Omega \in \{\mathbb{R}^d,\mathbb{T}^d\}$,} we write $\boldu\in C_{c,\sigma}^\infty(\Omega)$ if $\boldu \in C_c^\infty(\Omega)^d$ satisfies $\Div \boldu=0$ in $\Omega$. For $1\leq p<\infty$, we write $\Leb{p}_\sigma(\Omega)$ the closure of $C_{c,\sigma}^\infty(\Omega)$ under the $\Leb{p}$-norm.  For $1\leq p,q\leq \infty$, we write
\[ 
\norm{f}{\Leb{p}_t\Leb{q}_x((0,T)\times\Omega)}:=\norm{\norm{f(t)}{\Leb{q}_x(\Omega)}}{\Leb{p}_t(0,T)}.
\] 
We also define $\norm{f}{C^k_t H^s((0,T)\times\Omega)}$ for $k=0,1,\dots,$ and $s>0$.

\vspace{1ex}

%We drop the set in the parenthesis if it is ambient. 
 
For $f\in \Leb{1}(\Omega)$, we define the Fourier transform ($\Omega=\mathbb{R}^d$) and the Fourier coefficient ($\Omega=\mathbb{T}^d$) as follows: 
\[  
\begin{split}
\widehat{f}(\xi)&=\frac{1}{(2\pi)^{d/2}} \int_{\mathbb{R}^d} f(x) e^{-ix\cdot \xi} \myd{x},\quad \xi \in \mathbb{R}^d,\\
\widehat{f}(k)&=\frac{1}{(2\pi)^{d/2}} \int_{\mathbb{T}^d} f(x) e^{-ix\cdot k} \myd{x},\quad k \in \mathbb{Z}^d.
\end{split}
\]
Using this, we define the inhomogeneous Sobolev spaces $H^s(\Omega)$ of order $s$ with the following norms
\[   
\begin{split}
\norm{f}{H^s(\mathbb{R}^d)}^2 &= \int_{\mathbb{R}^d} (1+|\xi|^2)^s |\widehat{f}(\xi)|^2 d\xi<\infty,\\
\norm{f}{H^s(\mathbb{T}^d)}^2 &= \sum_{k\in \mathbb{Z}^d} (1+|k|^2)^s |\widehat{f}(k)|^2 <\infty.
\end{split}
\]
We write $\tSob{s}_\sigma(\Omega)$ the closure of $C^\infty_{c,\sigma}(\Omega)$ under $\tSob{s}$-norm.  The homogeneous Sobolev space $\thSob{s}(\Omega)$ consists of measurable functions $f$ for which
\[  
\begin{split}
\norm{f}{\thSob{s}(\mathbb{R}^d)}^2 &= \int_{\mathbb{R}^d} |\xi|^{2s} |\widehat{f}(\xi)|^2 d\xi<\infty, \\
\norm{f}{\thSob{s}(\mathbb{T}^d)}^2 &= \sum_{k\in\mathbb{Z}^d} |k|^{2s} |\widehat{f}(k)|^2<\infty.
\end{split}
\]

For $\boldu\in \mathcal{S}(\mathbb{R}^d)^d$, we define the Leray projection $\mathbb{P}$ by
\[ \widehat{\mathbb{P}(\boldu)}(\xi)= \left(I-\frac{\xi\otimes \xi}{|\xi|^2}\right)\widehat{\boldu}(\xi).\]
 It is well-known that the Leray projection is a bounded linear operator from $\Leb{q}(\mathbb{R}^d)^d$ into $\Leb{q}_\sigma(\mathbb{R}^d)$. Similarly, we can define the Leray projection for $\boldu\in C^\infty(\mathbb{T}^d)^d$ having {zero mean}. {The projection is a bounded linear operator from $\dot{L^q}(\mathbb{T}^d)^d$ to $\Leb{q}_\sigma(\mathbb{T}^d)$, where
 \[
\dot{L^q}(\mathbb{T}^d)^d= \left\{ \boldu \in \Leb{q}(\mathbb{T}^d)^d : \int_{\mathbb{T}^d} \boldu \myd{x}=0\right\}\quad \text{and}\quad  \dot{L^q_\sigma}(\mathbb{T}^d)= \left\{ \boldu \in \Leb{q}_\sigma(\mathbb{T}^d) : \int_{\mathbb{T}^d} \boldu \myd{x}=0\right\}
 \] (see \cite[Theorem 2.28]{RRS16}).}

%%%%%%%%%%%%%%%
\subsection{Main results}
%%%%%%%%%%%%%%%
We now describe the main results of this paper in detail. The first two results, Theorem \ref{thm:LWP} and Theorem \ref{thm:A}, concern local and global existence of strong solutions of \eqref{eq:MRE} on $\mathbb{T}^d$, $d=2,3$ for the non-resistive case $\eta=0$. In this case, Beekie-Friedlander-Vicol \cite{BFV22} prove{d} the local existence of $H^s$-solutions of \eqref{eq:MRE} when $\alpha \geq 0$ and $s>\frac{d}{2}+1$. In the same paper, they also give the continuation condition of strong solutions when $\alpha,s>\frac{d}{2}+1$. {Our first result is} a refinement of the local-in-time existence result of \cite{BFV22} by weakening the range of $s$ and the upper bound of solutions.

\begin{thm} \label{thm:LWP}
Let $d=2,3$. Suppose that either 
\eqn \label{eq:torus condition of LWP}
\begin{split}
\mathrm{(i)}\quad& \alpha>\frac{d}{2} \quad \text{and} \quad s\geq 1 \quad \text{or}\\
\mathrm{(ii)}\quad& 0\leq \alpha \leq \frac{d}{2} \quad \text{and} \quad s>\frac{d}{2}+1-\alpha
\end{split}
\een
holds and assume that $\boldb_{0}\in H^{s}_\sigma(\mathbb{T}^{d})$. Then there exists $T\geq (C\norm{\boldb_{0}}{H^{s}})^{-2}$ such that \eqref{eq:MRE} with $\eta=0$ has a unique solution $\boldb\in C([0,T);H^{s}(\mathbb{T}^{d}))$, with the associated zero-mean velocity $\boldu\in C([0,T);\dot{H}^{s+\alpha^*}(\mathbb{T}^{d}))\cap L^{2}((0,T);\dot{H}^{s+\alpha}(\mathbb{T}^{d}))$, where 
\[ \alpha^*=\begin{cases}
\alpha &\quad \text{if } (\alpha,s) \text{ satisfies (i)};\\
2\alpha-\frac{d}{2} &\quad \text{if } (\alpha,s) \text{ satisfies (ii)}.
\end{cases}
\] 
Moreover, $(\boldb,\boldu)$ satisfies
\begin{equation}\label{eq:torus-energy-ineq}
\sup_{\tau\in[0,t]}\norm{\boldb(\tau)}{L^{2}}^{2}+2\int^{t}_{0}\norm{\boldu(\tau)}{\dot{H}^{\alpha}}^{2}\myd{\tau}\leq \norm{\boldb_{0}}{L^{2}}^{2}.
\end{equation}
When (i) holds, $(\boldb,\boldu)$ also satisfies
\begin{equation}\label{eq:torus-exp}
\begin{aligned}
&\norm{\boldb(t)}{H^{s}}^{2}+\int^{t}_{0}\norm{\boldu(\tau)}{\dot{H}^{s+\alpha}}^{2}\myd{\tau} \\
&\leq \norm{\boldb_{0}}{H^{s}}^{2}\exp{\left(C\int^{t}_{0}\norm{\nabla\boldu(\tau)}{L^{\infty}}\myd{\tau}+t\norm{\boldb_{0}}{H^{1}}^{2}\exp{\left(C\int^{t}_{0}\norm{\nabla\boldu(\tau)}{L^{\infty}}\myd{\tau}\right)}\right)}
\end{aligned}
\end{equation}
for some constant $C=C(\alpha,s,d)>0$ and for all $t\in[0,T)$.
\end{thm}

\begin{rem*}
Compared to \cite[Theorem 2.2]{BFV22}, we reduce the range of $s$ from $s>\frac{d}{2}+1$ to \eqref{eq:torus condition of LWP} and the right-hand side of \eqref{eq:torus-exp} only contains $\norm{\nabla\boldu}{L^{\infty}}$, while the right-hand side {of (2.4) in} \cite[Theorem 2.2]{BFV22} contains $\norm{\nabla\boldu}{L^{\infty}}$ and $\norm{\nabla\boldb}{L^{\infty}}$.
\end{rem*}

{Our} next result concerns a sufficient condition for the continuation of solutions in terms of $\norm{\boldu}{L^{1}_{t}\dot{H}^{\frac{d}{2}+1}_{x}}$ when $\alpha>\frac{d}{2}$. As an application, we show that the local-in-time solution in Theorem \ref{thm:LWP} can be defined globally-in-time when $\alpha\geq \frac{d}{2}+1$. Therefore, we answer the question of the global existence part of \cite[Q1]{BFV22} when $\alpha\geq \frac{d}{2}+1$.

\begin{thm}\label{thm:A} Let $d=2,3$.
\begin{enumerate}[label=\textnormal{(\roman*)}]
\item Suppose that $\alpha>\frac{d}{2}$, $s\geq1$, and $\boldb_{0}\in H^{s}_\sigma(\mathbb{T}^{d})$. If $(\boldb,\boldu)$ is a local-in-time solution established in Theorem \ref{thm:LWP} and if
\begin{equation}\label{eq:BKM}
\int^{T}_{0}\norm{\boldu(t)}{\dot{H}^{\frac{d}{2}+1}(\mathbb{T}^d)}\myd{t}<\infty,
\end{equation}
then the solution can be extended beyond $t=T$.
\item If  $\alpha\geq\frac{d}{2}+1$ and $s\geq 1$, then the local-in-time solution in Theorem \ref{thm:LWP} is defined globally-in-time, i.e., $\boldb \in C([0,T];\tSob{s}(\mathbb{T}^d))$ for any $T>0$.
 \end{enumerate}
\end{thm}

The third result of this paper, Theorem \ref{thm:B}, concerns the existence, uniqueness, and the asymptotic behavior of mild solutions to \eqref{eq:MRE} in $\mathbb{R}^d$ for the resistive case $\eta>0$ for which we set $\eta=1$ for simplicity. We first find a proper space to deal with this case which usually comes from the scaling-invariant property of \eqref{eq:MRE}. One can see that a natural scaling corresponding to \eqref{eq:MRE} is 
\[
\boldb_\lambda (x,t)=\lambda^{\alpha+\beta-1} \boldb(\lambda x,\lambda^{2\beta} t),\quad \boldu_\lambda (x,t)=\lambda^{2\beta-1} \boldu(\lambda x, \lambda^{2\beta} t).
\]
Under this scaling property, we have
\[
\norm{\boldb_\lambda}{\Leb{p}}=\lambda^{\alpha+\beta-1-\frac{d}{p}} \norm{\boldb_0}{\Leb{p}}.
\]
From this, we classify the power of $\lambda$ into three cases:
\begin{enumerate}
\item subcritical : $\alpha+\beta-1-\frac{d}{p}>0$;
\item critical : $\alpha+\beta-1-\frac{d}{p}=0$;
\item supercritical : $\alpha+\beta-1-\frac{d}{p}<0$.
\end{enumerate} 

In this paper, we only deal with the critical case with the aim of showing solutions globally-in-time as in \cite{BBT12}. Formally, we can express a solution $\boldb$ of \eqref{eq:MRE} in the integral form:
\[
\boldb(t)=G_{\beta}(t)\boldb_0 -\int_0^t G_{\beta}(t-s) \Div(\boldu\otimes\boldb-\boldb\otimes \boldu)(s)\myd{s},
\]
where $\boldu=\Lambda^{-2\alpha} \Proj\Div(\boldb\otimes\boldb)$ and $G_{\beta}$ is {the fractional heat semigroup of order $\beta$ (see Subsection \ref{subsec:linear}).} Any solution satisfying this integral equation is called a mild solution.

{Now we are ready to present our third result of this paper.}
\begin{thm}\label{thm:B}
  Let $d \geq 2$ and let  $\alpha$, $\beta$ satisfy 
\[
\frac{1}{2}<\alpha< \frac{d+1}{2},\quad \beta>\frac{1}{2},\quad \alpha+\beta < d+1
\]
and 
\[
3\alpha+\beta<\frac{3d+4}{2}
\]
in addition if 
\[
\frac{d+2}{4}<\alpha<\frac{d+1}{2}.
\]
Let 
\[      
p=\frac{d}{\alpha+\beta-1}.
\]
Then for any $\boldb_0\in \Lebdiv{p}(\mathbb{R}^d)$, there exists $T_*>0$ such that there exists a unique {mild}  solution of \eqref{eq:MRE} satisfying $\boldb \in C([0,T_*];\Leb{p})\cap \mathcal{F}_{T_*}$, where
\[
\norm{\boldb}{\mathcal{F}_T} = \sup_{0<s\leq T} s^{\frac{d}{2\beta}\left(\frac{1}{p}-\frac{1}{q}\right)} \norm{\boldb(s)}{\Leb{q}},\quad \frac{1}{q}=\frac{1}{p}-\frac{2\beta-1}{3d}.
\] 
Moreover, if $\norm{\boldb_0}{\Leb{p}}$ is sufficiently small, then $T^{\ast}=\infty$ and 
\[
\lim_{t\rightarrow \infty}\norm{\boldb(t)}{\Leb{p}}=0.
\]
 \end{thm}

\begin{rem*}\leavevmode
\begin{enumerate}
\item[(i)] {In \cite{KK23}, the uniqueness of weak solutions was shown when $(\alpha,\beta)$ satisfies \eqref{eq:uniqueness-weak}. In particular, we need to assume that $\beta \geq 1$ to ensure the uniqueness of weak solutions. Theorem \ref{thm:B} shows that we can relax $\beta >\frac{1}{2}$ when we show the existence and uniqueness of mild solutions to \eqref{eq:MRE}.}
\item[(ii)] If we define $r$ by $\frac{1}{r}=\frac{2}{q}-\frac{2\alpha-1}{d}$, then $2<q<\infty$ and $1<r<\infty$. By the $\Leb{\frac{q}{2}}$-boundedness of the Leray projection and the Hardy-Littlewood-Sobolev inequality, we bound $\boldu$ as  
\begin{equation}\label{eq:boldu-control-p}
\norm{\boldu(t)}{\Leb{r}} =\norm{\Lambda^{-2\alpha} \Proj\Div(\boldb(t)\otimes\boldb(t))}{\Leb{q}} \apprle \norm{\boldb(t)}{\Leb{q}}^2.
\end{equation}
So, if $\norm{\boldb_0}{\Leb{p}}$ is sufficiently small, then \eqref{eq:boldu-control-p} and Theorem \ref{thm:B} imply that $\boldb\in \mathcal{F}_\infty$ and 
\[
   \norm{\boldu(t)}{\Leb{r}} \leq C\norm{\boldb(t)}{\Leb{q}}^2\leq C t^{-2\sigma},\quad \sigma=\frac{d}{2\beta}\left(\frac{1}{p}-\frac{1}{q}\right).
\]
On the other hand, suppose in addition that $\beta<\frac{d+1}{2}$. 
Then since $2<p<\infty$ and \eqref{eq:boldu-control-p} holds for $(s,p)$ instead of $(r,q)$, where $\frac{1}{s}=\frac{2}{p}-\frac{2\alpha-1}{d}$, we have
\[ 
\lim_{t\rightarrow\infty} \left(\norm{\boldu(t)}{\Leb{s}}+\norm{\boldb(t)}{\Leb{p}}\right)=0.
\]

\end{enumerate}
\end{rem*}

%%%%%%%%%%%%%%
\section{Preliminaries}\label{sec:prelim}
%%%%%%%%%%%%%%%
This section consists of four parts. We first give some elementary inequalities in Subsection \ref{subsec:inequalities}. We next  define fractional Laplacians on $\mathbb{R}^{d}$ and $\mathbb{T}^{d}$ in Subsection \ref{subsec:Fractional Laplacian}. We also recall Riesz potentials and Hardy-Littlewood-Sobolev inequality. After that, we proceed to give {several estimates in} Sobolev spaces on $\mathbb{T}^{d}$ in Subsection \ref{subsec:Sobolev}. In Subsection \ref{subsec:linear}, we list some bounds of linear semigroups for fractional heat operators on $\mathbb{R}^{d}$  that {are} used in the construction of solutions to \eqref{eq:MRE} when $\eta>0$.

%%%%%%%%%%%%%%%%%%%%%%%%%%%%%%%
\subsection{Some inequalities} \label{subsec:inequalities}
%%%%%%%%%%%%%%%%%%%%%%%%%%%%%%
We recall a few inequalities which will be {used} repeatedly when we prove our results. However, we will not refer them whenever it is obvious to use them. 
\begin{enumerate}[]
\item Young's inequality: if $p>1$ and $q>1$ satisfy $\frac{1}{p}+\frac{1}{q}=1$, then
\[
ab\leq \frac{a^p}{p}+\frac{b^q}{q} \quad \text{for $a\ge 0$ and $b\ge 0$.}
\]
\item H\"older's inequality: if $f\in L^{p}$ and $g\in L^{q}$, then $fg\in L^{1}$ with 
\[
\|fg\|_{L^{1}}\leq \|f\|_{L^{p}}\|g\|_{L^{q}}, \quad \frac{1}{p}+\frac{1}{q}=1.
\]
\item Gr\"onwall's inequality \cite[Page 624]{Evans}: Let $\eta$ be a nonnegative, absolutely continuous function on $[0,T]$ satisfying for a.e. $t$ the differential inequality
\[
\eta'(t) \leq \phi(t)\eta(t) +\psi(t),
\]
where $\phi(t)$ and $\psi(t)$ are nonnegative  and integrable functions on $[0,T]$. Then, 
\[
\eta(t) \leq \left(\eta(0)+\int^{t}_{0}\psi(s)ds\right)\exp\left[\int^{t}_{0}\phi(s)ds\right].
\]
\item \cite[Lemma 3]{BBT12}: for $0<a,b<1$,
\eqn \label{eq:beta-function-estimate}
\int_0^t (t-s)^{-a} s^{-b} ds \leq C t^{1-a-b}.
\een
\end{enumerate}

%%%%%%%%%%%%%%%%%%
\subsection{Fractional Laplacian and Riesz potentials}\label{subsec:Fractional Laplacian}
When the domain is $\mathbb{R}^d$, the fractional Laplacian $\Lambda^\gamma = (\sqrt{-\Delta})^{\gamma}$, $\gamma>0$ is defined by the Fourier transform
\begin{equation}\label{FL1}
\widehat{\Lambda^\gamma f}(\xi)=|\xi|^\gamma \widehat{f}(\xi),\quad \xi \in \mathbb{R}^d.
\end{equation}
Similarly, when the domain is $\mathbb{T}^d$, the fractional Laplacian $\Lambda^\gamma$ is defined as
\begin{equation}\label{FL2}
\widehat{\Lambda^\gamma f}(k)=|k|^\gamma \widehat{f}(k),\quad k\in\mathbb{Z}^d.
\end{equation}

For $0<\gamma<d$, we define the Riesz potentials of order $\gamma$ by
\[  
\Lambda^{-\gamma} f(x)=C_{\gamma}\int_{\mathbb{R}^d} \frac{f(y)}{|x-y|^{d-\gamma}}dy,\quad x\in\mathbb{R}^d,
\]
where $C_{\gamma}$ is a normalizing constant. Using this, we state the Hardy-Littlewood-Sobolev inequality \cite[p.106]{Lieb01}.

\begin{prop}\label{thm:HLS}
Let $0<\gamma<d$ and let $1<p<q<\infty$ satisfy
\[ 
\frac{1}{q}=\frac{1}{p}-\frac{\gamma}{d}.
\]
Then there exists a constant $C=C(d,\gamma,p)>0$ such that 
\[ 
\norm{\Lambda^{-\gamma} f}{\Leb{q}(\mathbb{R}^d)}\leq C \norm{f}{\Leb{p}(\mathbb{R}^d)}
\]
for all $f\in \mathcal{S}(\mathbb{R}^d)$ and so the operator $\Lambda^{-\gamma}$ uniquely extended to a bounded linear operator from $\Leb{p}(\mathbb{R}^d)$ to $\Leb{q}(\mathbb{R}^d)$.
\end{prop}

%%%%%%%%%%%%%%%%%%%%%%%%%%%%%%%%%%%%%%
\subsection{Sobolev spaces on $\mathbb{T}^{d}$}\label{subsec:Sobolev}
%%%%%%%%%%%%%%%%%%%%%%%%%%%%%%%%%%%%%%
In this subsection, we recall embedding and density results of Sobolev spaces $\tSob{s}(\mathbb{T}^d)$, fractional Leibniz rules, Kato-Ponce commutator estimates, and prove logarithmic type Sobolev inequality. We first recall an embedding theorem and density results on Sobolev spaces.

\begin{prop}\label{prop:density-embedding}\leavevmode
\begin{enumerate}
\item[\rm (i)] If $s>\frac{d}{2}$, then $\tSob{s}(\mathbb{T}^d)\hookrightarrow \Leb{\infty}(\mathbb{T}^d)$.
\item[\rm (ii)] If $0\leq s_1\leq s_2$ and $1<q<\infty$ satisfy
\[
 \frac{1}{2}-\frac{s_2}{d}\leq \frac{1}{q}-\frac{s_1}{d},
\]
then 
\[  
\norm{u}{\Leb{q}}+\norm{\Lambda^{s_1} u}{\Leb{q}}\apprle \norm{u}{\tSob{s_2}}.
\]
\item[\rm (iii)] If $f\in H^s(\mathbb{T}^d)$ has {zero mean} for $s>0$, then 
\[
\norm{f}{\Leb{2}}\leq \norm{f}{\dot{H}^{s}}\approx \norm{f}{\tSob{s}}.
\]
\noindent
\item[\rm (iv)] If $s>0$ and $f\in \tSob{s}(\mathbb{T}^d)$ has {zero mean}, then there exists a sequence $\{f_k\} \subset C^\infty(\mathbb{T}^d)$ with $\displaystyle \int_{\mathbb{T}^d} f_k \myd{x}=0$ such that $f_k\rightarrow f$ in $\tSob{s}(\mathbb{T}^d)$. 
\end{enumerate}
\end{prop}

\begin{proof}
(i) and (ii):  These are well-known Sobolev embeddings,  see B\'{e}nyi-Oh \cite{BO13} or Cirant-Goffi \cite[Lemma 2.5]{CG19}  for the proof.

\vspace{1ex}

\noindent
(iii) For a zero-mean vector field $f$, since $\widehat{f}(0)=0$, it follows from Parseval's identity that 
\begin{align*}
\norm{f}{\Leb{2}}^{2}=\sum_{k\in\mathbb{Z}^{d}\setminus \{0\}}\left|\widehat{f}(k)\right|^{2} &\leq \sum_{k\in\mathbb{Z}^{d}}|k|^{2s}\left|\widehat{f}(k)\right|^{2}=:\norm{f}{\dot{H}^{s}}^{2} \\
&\approx \sum_{k\in\mathbb{Z}^{d}}(1+|k|^{2})^{s}\left|\widehat{f}(k)\right|^{2}=\norm{f}{H^{s}}^{2}.
\end{align*}
\noindent
(iv) Since $C^\infty(\mathbb{T}^d)$ is dense in $\tSob{s}(\mathbb{T}^d)$, there exists a sequence $\phi_k \in C^\infty(\mathbb{T}^d)$ such that $\norm{\phi_k-f}{\tSob{s}}\rightarrow 0$ as $k\rightarrow \infty$. Define 
\[
f_k = \phi_k - \frac{1}{|\mathbb{T}^d|}\int_{\mathbb{T}^d} \phi_k \myd{x}.
\]
Then $\displaystyle \int_{\mathbb{T}^d} f_k \myd{x}=0$ and so $f_k-f$ has {zero mean}. Hence we deduce from (iii) that 
\[
\norm{f_k-f}{\tSob{s}}\leq C\norm{f_k-f}{\dot{H}^s} \leq C\norm{\phi_k-f}{\dot{H}^s} \leq C \norm{\phi_k-f}{\tSob{s}}\rightarrow 0
\]
as $k\rightarrow \infty$. This completes the proof of Proposition \ref{prop:density-embedding}.
\end{proof}

We recall the following fractional Leibniz rules and Kato-Ponce commutator estimates (see e.g. \cite{KP88,KPV91,Ju04,L19,CK19} and references therein).

\begin{prop} Let $s>0$, $1<p<\infty$, and $1<p_{1},p_{2},p_{3},p_{4}\leq\infty$ satisfy $\frac{1}{p}=\frac{1}{p_{1}}+\frac{1}{p_{2}}=\frac{1}{p_{3}}+\frac{1}{p_{4}}$. Then we have
\begin{align}\label{eq:product-estimate}
\norm{\Lambda^{s}(fg)}{L^{p}}&\apprle \left(\norm{\Lambda^{s}f}{L^{p_{1}}}\norm{g}{L^{p_{2}}}+\norm{f}{L^{p_{3}}}\norm{\Lambda^{s}g}{L^{p_{4}}}\right),\\
\label{eq:commutator-estimate}
\norm{\Lambda^{s}(fg)-f(\Lambda^{s}g)}{L^{p}}&\apprle  \left(\norm{\Lambda^{s}f}{L^{p_{1}}}\norm{g}{L^{p_{2}}}+\norm{\nabla f}{L^{p_{3}}}\norm{\Lambda^{s-1}g}{L^{p_{4}}}\right),
\end{align}
for all $f,g \in C^\infty(\mathbb{T}^d)$.
\end{prop} 

The following logarithmic type Sobolev inequality is usually proved in $\mathbb{R}^{d}$ as one can see, for example, in \cite{BKM84,BG80,BW80,KT00,KOT02,C02,O03}. However, we are not able to find a proper reference of showing a similar inequality in $\mathbb{T}^{d}$. Hence we give a proof for the sake of convenience.

\begin{prop}\label{prop:BKM-estimate}
Let $s>\frac{d}{2}$. Then we have
\[
\norm{f}{\Leb{\infty}}\apprle 1+\norm{f}{\thSob{\frac{d}{2}}}\log(e+\norm{f}{\thSob{s}})
\]
for all $f\in \dot{H}^s$ {having zero mean}.
\end{prop}

\begin{proof}
By Proposition \ref{prop:density-embedding} (iv), we may assume that $f$ is smooth {having zero mean}. Using the Fourier series expansion and Cauchy-Schwarz inequality, we have 
\begin{align*}
\left|f(x)\right|&\leq\sum_{k\in\mathbb{Z}^{d}\setminus\{0\}}\left|\widehat{f}(k)\right| = \sum_{1\leq|k|\leq R}|k|^{\frac{d}{2}}\left|\widehat{f}(k)\right|\cdot |k|^{-\frac{d}{2}}+\sum_{|k|>R}|k|^{s}\left|\widehat{f}(k)\right|\cdot|k|^{-s} \\
&\lesssim \norm{f}{\dot{H}^{\frac{d}{2}}}\left(\sum_{1\leq|k|\leq R}|k|^{-d}\right)^{\frac{1}{2}}+\norm{f}{\dot{H}^{s}}\left(\sum_{|k|>R}|k|^{-2s}\right)^{\frac{1}{2}},
\end{align*}
where $R$ is any real number greater than 1. To estimate the first summation term on the right-hand side, we observe that 
\[
\sum_{1\leq|k|\leq R}|k|^{-d}\leq \sum_{1\leq r\leq R}r^{-d}|A_{r}|,\quad\text{where } A_{r}=\left\{k\in\mathbb{Z}^{d}: r^{2}\leq |k|^{2}<(r+1)^{2}\right\}
\]
and we denote $|A|$ as the number of elements of a set $A$. Since the number of integers in each shell approximates its area (or volume), we have $|A_{r}|\lesssim r^{d-1}$. So we derive the following inequality: 
\begin{equation}\label{eq:estimate-1}
\sum_{1\leq|k|\leq R}|k|^{-d} \lesssim \sum_{1\leq r\leq R}r^{-d}\cdot r^{d-1}\lesssim 1+\int^{R}_{1}x^{-1}dx\lesssim 1+\log{R}.
\end{equation}
Since $s>\frac{d}{2}$, we estimate the second summation by
\begin{equation}\label{eq:estimate-2}
\sum_{|k|>R}|k|^{-2s}\lesssim \sum_{r>R}r^{-2s}\cdot r^{d-1}\lesssim   \int^{\infty}_{R}x^{d-1-2s}dx\lesssim R^{d-2s}.
\end{equation} 
Hence by \eqref{eq:estimate-1} and \eqref{eq:estimate-2}, we obtain
\[
\norm{f}{L^{\infty}}\lesssim \norm{f}{\dot{H}^{\frac{d}{2}}}\left(1+\log{R}\right)+\norm{f}{\dot{H}^{s}}R^{\frac{d}{2}-s}.
\]
By choosing $R\ge1$ so that $R^{s-\frac{d}{2}}=\norm{f}{\dot{H}^{s}}/\norm{f}{\dot{H}^{\frac{d}{2}}}$, we have
\eqn \label{eq:beforeBKM}
\norm{f}{\Leb{\infty}}\apprle \norm{f}{\thSob{\frac{d}{2}}}\log \left(e+\norm{f}{\dot{H}^{s}}/\norm{f}{\dot{H}^{\frac{d}{2}}}\right).
\een
Since 
\[
      x \left(\log (e+ y/x)\right)\leq \begin{cases}
       1+x\log(e+y) &\text{if } 0<x\leq 1 \\
       x \log (e+y)&\text{if } x\geq 1
       \end{cases}
\]
holds for all $x,y>0$, {it follows from \eqref{eq:beforeBKM} that}
\begin{align*}
\norm{f}{\Leb{\infty}}&\apprle 1+\norm{f}{\thSob{\frac{d}{2}}}\log(e+\norm{f}{\thSob{s}}).
\end{align*}
This completes the proof of Proposition \ref{prop:BKM-estimate}. 
\end{proof}

%%%%%%%%%%%%%%%%%%%%%%%%%%%%%%%%%%%%%%%%%%
\subsection{Fractional heat equations on $\mathbb{R}^d$}\label{subsec:linear}
%%%%%%%%%%%%%%%%%%%%%%%%%%%%%%%%%%%%%%%%%%
In this subsection, we summarize several known results on derivative estimates for fractional heat semigroups on $\mathbb{R}^d$ and their $\Leb{p}$-mapping properties.

For $\beta>0$, we consider the following fractional heat equations in $(0,\infty)\times \mathbb{R}^d$: 
\begin{equation} \label{eq:linear fractional heal eq}
\left\{
\begin{aligned}
\partial_t \boldb+\Lambda^{2\beta}\boldb&=0&&\quad\text{in } (0,\infty)\times\mathbb{R}^d,\\
\boldb&=\boldb_0&&\quad\text{on } \{t=0\}\times \mathbb{R}^d.
\end{aligned}
\right.
\end{equation}
Then, $\boldb(t,x)$ is expressed as
\[   
\boldb(t,x)=\mathcal{F}^{-1}\left(e^{-t|\xi|^{2\beta}} \mathcal{F}{\boldb_0}(\xi)\right)=(K_\beta (t,\cdot) * \boldb_0)(x), \quad K_\beta(t,x)=\mathcal{F}^{-1}\left(e^{-t|\xi|^{2\beta}}\right)(x).
\]
Let us write $G_\beta (t)$ the convolution operator with the kernel $K_\beta (t,x)$, i.e., $G_\beta (t)f=K_\beta (t,\cdot)*f$.

\begin{prop}[{\cite[Lemma 3.1]{MYZ08}}]\label{prop:linear-estimate}
Let $1\leq p\leq q\leq \infty$ and $\beta,\gamma>0$. Then 
\begin{align}
\norm{G_\beta(t)\varphi}{\Leb{q}}&\apprle t^{-\frac{d}{2\beta}\left(\frac{1}{p}-\frac{1}{q}\right)}\norm{\varphi}{\Leb{p}},\\
\norm{\Lambda^\gamma G_\beta(t)\varphi}{\Leb{q}}&\apprle t^{-\frac{\gamma}{2\beta}-\frac{d}{2\beta}\left(\frac{1}{p}-\frac{1}{q}\right)}\norm{\varphi}{\Leb{p}}
\end{align}
for all $\varphi \in \Leb{p}(\mathbb{R}^d)$.
\end{prop}

Proposition \ref{prop:linear-estimate} gives the following limit behavior of  $G_\beta(t)\boldb_0$ near $t=0$.

\begin{prop}\label{prop:time-control}
Let $\beta>0$ and $1\leq p< q<\infty$. Then
\[
 \lim_{t\rightarrow 0+} t^{\frac{d}{2\beta}\left(\frac{1}{p}-\frac{1}{q} \right)} \norm{G_\beta(t)\boldb_0}{\Leb{q}}=0.
\]
\end{prop}

\begin{proof}
Although this can be proven using the argument in \cite{CF08}, we give a detail of the proof for the sake of completeness.  First of all, by Proposition \ref{prop:linear-estimate}, we have  
\[
 t^{\frac{d}{2\beta}\left(\frac{1}{p}-\frac{1}{q} \right)} \norm{G_\beta(t)\boldb_0}{\Leb{q}}\apprle \norm{\boldb_0}{\Leb{p}}
\]
for $t>0$. We now choose $\boldb_0^k \in C_c^\infty$ so that $\norm{\boldb_0^k-\boldb_0}{\Leb{p}}\rightarrow 0$ as $k\rightarrow\infty$. For this choice $\boldb_0^k$, we have 
\[
  \norm{G_\beta(t)\boldb_0^k}{\Leb{q}}\apprle \norm{\boldb_0^k}{\Leb{q}}.
\]
Hence 
\[
\lim_{t\rightarrow 0+}t^{\frac{d}{2\beta}\left(\frac{1}{p}-\frac{1}{q} \right)} \norm{G_\beta(t)\boldb_0^k}{\Leb{q}}=0
\]
for each $k$. Moreover, $t^{\frac{d}{2\beta}\left(\frac{1}{p}-\frac{1}{q} \right)}G_\beta(t)\boldb_0^k$ is in $BC(0,\infty;\Leb{q})$ and it converges to a function $t^{\frac{d}{2\beta}\left(\frac{1}{p}-\frac{1}{q} \right)}G_\beta(t)\boldb_0$ uniformly in any interval $(0,T)$ for all $T>0$. This completes the proof of Proposition \ref{prop:time-control}. 
\end{proof}

We finally recall the semigroup property of the fractional heat operator.

\begin{prop}\label{prop:semigroup}
For $t_1,t_2>0$, $1\leq p< \infty$, and $\beta>0$, we have 
\[  G_{\beta}(t_1+t_2)\varphi = G_{\beta}(t_1)(G_{\beta}(t_2)\varphi)\]
for all $\varphi \in \Leb{p}(\mathbb{R}^d)$.
\end{prop}

%%%%%%%%%%%%%%%%%%%%%%%%%%%%%%%%%%%%%%%%%%%%%%%
\section{Non-resistive case: Proof of Theorem \ref{thm:LWP} and \ref{thm:A}}\label{sec:4}
%%%%%%%%%%%%%%%%%%%%%%%%%%%%%%%%%%%%%%%%%%%%%%%

\subsection{Proof of Theorem \ref{thm:LWP}}
We first recall \eqref{eq:MRE} with $\eta=0$ {and $\nu=1$}:
\begin{equation}\label{eq:MRE N}
\left\{
\begin{alignedat}{2}
\Lambda^{2\alpha} \boldu +\nabla \pbar &=(\boldb\cdot \nabla)\boldb&&\quad \text{in } \mathbb{T}^d\times (0,T),\\
\partial_t \boldb+(\boldu\cdot \nabla)\boldb &=(\boldb\cdot \nabla)\boldu &&\quad \text{in }\mathbb{T}^d \times(0,T),\\
\Div \boldu =\Div \boldb &=0 &&\quad \text{in } \mathbb{T}^d\times (0,T),\\
\boldb(\cdot,0)&=\boldb_0 &&\quad \text{on } \mathbb{T}^d.
\end{alignedat}\right.
\end{equation}

In order to prove Theorem \ref{thm:LWP}, it {suffices} to obtain a priori estimates for strong solutions because {the existence of solutions follows by the standard approximation argument, and} the uniqueness of solutions $(\boldb, \boldu)$ in $L^{\infty}([0,T);L^{2}(\mathbb{T}^{d}))\times L^{2}([0,T);\dot{H}^{\alpha}(\mathbb{T}^{d}))$ follows by a similar process to obtain the a priori estimates.

By multiplying  the first equation in \eqref{eq:MRE N} by $\boldu$ and the second equation by $\boldb$, by integrating the resulting equations on $\mathbb{T}^d$, and using $\Div \boldu=\Div\boldb=0$ in $\mathbb{T}^d$, we deduce
\begin{equation}\label{eq:torus-L2}
\frac{1}{2}\frac{d}{dt}\norm{\boldb}{\Leb{2}}^{2}+\norm{\boldu}{\dot{H}^{\alpha}}^{2}=0
\end{equation} 
from which we arrive at \eqref{eq:torus-energy-ineq}. 

By multiplying $-\Delta \boldu$ and $-\Delta \boldb$ to the first equation and the second equation in \eqref{eq:MRE N}, respectively, we also have 
\begin{equation}\label{eq:torus-H1}
\begin{aligned}
\frac{1}{2}\frac{d}{dt}\norm{\boldb}{\dot{H}^{1}}^{2}+\norm{\boldu}{\dot{H}^{1+\alpha}}^{2}&=\int_{\mathbb{T}^d}\Delta\boldb\cdot (\boldu\cdot\nabla\boldb)\myd{x}-\int_{\mathbb{T}^d}\Delta\boldb\cdot (\boldb\cdot\nabla\boldu)\myd{x}-\int_{\mathbb{T}^d}\Delta\boldu\cdot (\boldb\cdot\nabla\boldb)\myd{x} \\
&=-\int_{\mathbb{T}^d}\partial_{k}b^{j}\partial_{k}u^{i}\partial_{i}b^{j}+\int_{\mathbb{T}^d}\partial_{k}b^{j}\partial_{k}b^{i}\partial_{i}u^{j}+\int_{\mathbb{T}^d}\partial_{k}u^{j}\partial_{k}b^{i}\partial_{i}b^{j} \\
&\leq C\norm{\nabla\boldu}{L^{\infty}}\norm{\boldb}{\dot{H}^{1}}^{2} ,
\end{aligned}
\end{equation}
where $\Div\boldu=\Div\boldb=0$ are used in the second equality. By Gr\"onwall's inequality with \eqref{eq:torus-L2} and \eqref{eq:torus-H1}, we have
\begin{equation}\label{eq:torus-H1-exp}
\norm{\boldb(t)}{{H}^{1}}^{2}+2\int^{t}_{0}\norm{\boldu(\tau)}{\dot{H}^{1+\alpha}}^{2}\myd{\tau}\leq \norm{\boldb_{0}}{H^{1}}^{2}\exp{\left(C\int^{t}_{0}\norm{\nabla \boldu(\tau)}{L^{\infty}}\myd{\tau}\right)}.
\end{equation}

To obtain $\dot{H}^{s}$ energy estimates, we write
\begin{equation}\label{eq:torus-Hs-pre}
\begin{aligned}
\frac{1}{2}\frac{d}{dt}\norm{\boldb}{\dot{H}^{s}}^{2}+\norm{\boldu}{\dot{H}^{s+\alpha}}^{2}&=-\int_{\mathbb{T}^d}\Lambda^{s}\boldb\cdot\left[\Lambda^{s},\boldu\cdot\nabla\right]\boldb\myd{x}+\int_{\mathbb{T}^d}\Lambda^{s}\boldb\cdot\left[\Lambda^{s},\boldb\cdot\nabla\right]\boldu \myd{x} \\
&+\int_{\mathbb{T}^d}\Lambda^{s}\boldu\cdot\left[\Lambda^{s},\boldb\cdot\nabla\right]\boldb \myd{x},
\end{aligned}
\end{equation}
where $[A,B]=AB-BA$ is the usual commutator notation. We now recall two conditions of $(\alpha, s)$ specified in Theorem \ref{thm:LWP}: either
\eqn \label{LWP Condition 1}
\alpha>\frac{d}{2} \quad \text{and} \quad s\geq 1
\een
or 
\eqn \label{LWP Condition 2}
0\leq \alpha \leq \frac{d}{2} \quad \text{and} \quad s>\frac{d}{2}+1-\alpha.
\een

{\noindent \emph{Case 1}.} Suppose first that (\ref{LWP Condition 1}) holds. By H\"older's inequality, Kato-Ponce commutator estimate \eqref{eq:commutator-estimate}, and Proposition \ref{prop:density-embedding} (i), we bound the first two terms on the right-hand side of (\ref{eq:torus-Hs-pre}): 
\[
\begin{aligned}
&-\int_{\mathbb{T}^d}\Lambda^{s}\boldb\cdot\left[\Lambda^{s},\boldu\cdot\nabla\right]\boldb\myd{x}+\int_{\mathbb{T}^d}\Lambda^{s}\boldb\cdot\left[\Lambda^{s},\boldb\cdot\nabla\right]\boldu \myd{x}\\
&\apprle \norm{\boldb}{\dot{H}^{s}}\left(\norm{\left[\Lambda^{s},\boldu\cdot\nabla\right]\boldb}{L^{2}}+\norm{\left[\Lambda^{s},\boldb\cdot\nabla\right]\boldu}{L^{2}}\right) \\
&\apprle \norm{\boldb}{\dot{H}^{s}}\left(\norm{\Lambda^{s}\boldu}{L^{\infty}}\norm{\nabla\boldb}{L^{2}}+\norm{\nabla\boldu}{L^{\infty}}\norm{\Lambda^{s}\boldb}{L^{2}}\right) \\
&\apprle \norm{\boldu}{\dot{H}^{s+\alpha}}\norm{\boldb}{\dot{H}^{1}}\norm{\boldb}{\dot{H}^{s}}+\norm{\nabla\boldu}{L^{\infty}}\norm{\boldb}{\dot{H}^{s}}^{2}.
\end{aligned}
\]
Using the divergence-free condition on $\boldb$, H\"older's inequality,  fractional Leibniz rule \eqref{eq:product-estimate}, we estimate the last term on the right-hand side of (\ref{eq:torus-Hs-pre}) as
\[
\begin{aligned}
\int_{\mathbb{T}^d}\Lambda^{s}\boldu\cdot\left[\Lambda^{s},\boldb\cdot\nabla\right]\boldb \myd{x}&=-\int _{\mathbb{T}^d}\partial_i\Lambda^{s}u^j  (\Lambda^s(b^ib^j)-b^i \Lambda^s b^j ) \myd{x} \\
&\apprle \norm{\nabla\Lambda^{s}\boldu}{L^{\frac{d}{1-\epsilon}}}\left(\norm{\Lambda^s (\boldb\otimes\boldb)}{\Leb{\frac{d}{d-1+\varepsilon}}} +\norm{\boldb \otimes \Lambda^s \boldb}{\Leb{\frac{d}{d-1+\varepsilon}}}\right)\\
&\apprle \norm{\nabla\Lambda^{s}\boldu}{L^{\frac{d}{1-\epsilon}}}\norm{\Lambda^{s}\boldb}{L^{2}}\norm{\boldb}{L^{\frac{d}{d/2-1+\epsilon}}}.
\end{aligned}
\]
Since 
\eqn \label{eq:torus-Sobolev em}
\norm{\nabla\Lambda^{s}\boldu}{L^{\frac{d}{1-\epsilon}}} \apprle \norm{\nabla\Lambda^{s}\boldu}{\dot{H}^{\frac{d}{2}-1+\epsilon}} \leq \norm{\boldu}{\dot{H}^{\frac{d}{2}+\epsilon+s}},\quad \norm{\boldb}{L^{\frac{d}{d/2-1+\epsilon}}} \apprle \norm{\boldb}{\dot{H}^{1-\epsilon}},
\een
we choose $0<\varepsilon<\min\{1,\alpha-d/2\}$ to derive  
\[
\int_{\mathbb{T}^d}\Lambda^{s}\boldu\cdot\left[\Lambda^{s},\boldb\cdot\nabla\right]\boldb \myd{x} \apprle \norm{\boldu}{\dot{H}^{s+\alpha}} \norm{\boldb}{H^{1}}  \norm{\boldb}{\dot{H}^{s}}. 
\]
Hence we obtain
\begin{equation}\label{eq:torus-Hs}
\frac{1}{2}\frac{d}{dt}\norm{\boldb}{\dot{H}^{s}}^{2}+\norm{\boldu}{\dot{H}^{s+\alpha}}^{2}\apprle \norm{\nabla\boldu}{L^{\infty}}\norm{\boldb}{\dot{H}^{s}}^{2}+\norm{\boldu}{\dot{H}^{s+\alpha}}\norm{\boldb}{H^{1}}\norm{\boldb}{\dot{H}^{s}}.
\end{equation}
Since $\norm{\nabla \boldu}{\Leb{\infty}}\apprle \norm{\boldu}{\thSob{s+\alpha}}$ and $s\ge 1$, \eqref{eq:torus-L2} and \eqref{eq:torus-Hs} with Young's inequality imply 
\begin{equation}\label{eq:local-existence-estimate}
\frac{d}{dt}\norm{\boldb}{H^{s}}^{2}+\norm{\boldu}{\dot{H}^{s+\alpha}}^{2}\apprle \norm{\boldb}{H^{s}}^{4}.
\end{equation}
Therefore, \eqref{eq:local-existence-estimate} implies the local existence of solutions $(\boldb,\boldu) \in C([0,T);H^{s}(\mathbb{T}^{d}))\times L^{2}([0,T);\dot{H}^{s+\alpha}(\mathbb{T}^{d}))$.  Moreover, Gr\"onwall's inequality with \eqref{eq:torus-L2} and \eqref{eq:torus-Hs} yields that
\begin{equation}\label{eq:torus-Hs-exp}
\norm{\boldb(t)}{H^{s}}^{2}+\int^{t}_{0}\norm{\boldu(\tau)}{\dot{H}^{s+\alpha}}^{2}\myd{\tau}\leq \norm{\boldb_{0}}{H^{s}}^{2}\exp{\left(C\int^{t}_{0}\norm{\nabla \boldu(\tau)}{L^{\infty}}+\norm{\boldb(\tau)}{H^{1}}^{2}\myd{\tau}\right)}.
\end{equation}
Hence the desired estimate \eqref{eq:torus-exp} follows from \eqref{eq:torus-H1-exp} and \eqref{eq:torus-Hs-exp}.

Finally, we verify that $\boldu\in C([0,T);\dot{H}^{s+\alpha}(\mathbb{T}^{d}))$. By the $\Leb{q}$-boundedness of Leray projection, and the fractional Leibniz rule \eqref{eq:product-estimate} with $0<\epsilon<\min\{1,\alpha-\frac{d}{2}\}$, we have
\begin{equation}\label{eq:u-estimate-b}
\begin{aligned}
\norm{\Lambda^{s+\alpha}\boldu(t)}{L^{2}}&=\norm{\Lambda^{s-\alpha}\mathbb{P}\Div(\boldb\otimes\boldb)(t)}{L^{2}}\\
&\apprle \norm{\Lambda^{s+1-\frac{d}{2}-\epsilon}\left(\boldb(t)\otimes\boldb(t)\right)}{L^{2}}\\
&\apprle \norm{\Lambda^{s+1-\frac{d}{2}-\epsilon}\boldb(t)}{L^{\frac{d}{1-\epsilon}}}\norm{\boldb(t)}{L^{\frac{d}{d/2-1+\epsilon}}} \apprle \norm{\boldb(t)}{H^{s}}^{2},
\end{aligned}
\end{equation}
where we use (\ref{LWP Condition 1}) and the same embedding relations in \eqref{eq:torus-Sobolev em}. {It follows from a similar process to obtain \eqref{eq:u-estimate-b} and $\boldb\in C([0,T);H^{s}(\mathbb{T}^{d}))$ that 
\[
\norm{\Lambda^{s+\alpha}\boldu(t+h)-\Lambda^{s+\alpha}\boldu(t)}{L^{2}}\apprle \left(\norm{\boldb(t+h)}{H^{s}}+\norm{\boldb(t)}{H^{s}}\right)\norm{\boldb(t+h)-\boldb(t)}{H^{s}}\longrightarrow 0
\]
as $h$ tends to zero.} This completes the proof of Theorem \ref{thm:A} when $(\alpha,s)$ satisfies (\ref{LWP Condition 1}).\medskip

{\noindent \emph{Case 2}.}  We now consider the second case (\ref{LWP Condition 2}). {Let $2\leq p<\infty$ and $2<q\leq\infty$ such that
\begin{equation}\label{eq:p-q}
q=\begin{cases}
\infty \quad &\text{if}\ \ \alpha=0,\\
\frac{d}{\alpha} \quad &\text{if}\ \ \alpha\in(0,\frac{d}{2}),\\
2+\epsilon \quad &\text{if}\ \ \alpha=\frac{d}{2},
\end{cases}\qquad \frac{1}{p}=\frac{1}{2}-\frac{1}{q},
\end{equation}
where $\epsilon>0$ is a sufficiently small constant such that $\frac{1}{2+\epsilon}>\frac{1}{2}-\frac{s-1}{d}$.} Since $\boldu$ has {zero mean}, it follows from Proposition \ref{prop:density-embedding} {(i) and} (ii) that 
\begin{equation}\label{eq:joint-embedding}  \norm{\Lambda^s \boldu}{\Leb{p}}\apprle \norm{\boldu}{\tSob{s+\alpha}}\approx \norm{\boldu}{\thSob{s+\alpha}}\quad\text{and}\quad \norm{\nabla \boldb}{\Leb{q}}\apprle \norm{\boldb}{\tSob{s}}.
\end{equation}
By H\"older's inequality, Kato-Ponce commutator estimate \eqref{eq:commutator-estimate}, and \eqref{eq:joint-embedding}, we derive
\begin{equation}\label{eq:torus-Hs-cor}
\begin{aligned}
\frac{1}{2}\frac{d}{dt}\norm{\boldb}{\thSob{s}}^2+\norm{\boldu}{\thSob{s+\alpha}}^2&\apprle \norm{\nabla \boldu}{\Leb{\infty}}\norm{\boldb}{\thSob{s}}^2+\norm{\Lambda^s\boldu}{\Leb{p}}\norm{\nabla \boldb}{\Leb{q}}\norm{\boldb}{\thSob{s}}\\
&\apprle\norm{\boldu}{\thSob{s+\alpha}}\norm{\boldb}{\tSob{s}}^2,
\end{aligned}
\end{equation}
where we use $\norm{\nabla \boldu}{\Leb{\infty}}\apprle \norm{\boldu}{\thSob{s+\alpha}}$. By Young's inequality with \eqref{eq:torus-L2}, we arrive at the same bound \eqref{eq:local-existence-estimate}, which implies the local existence of a unique solution. By using Gr\"onwall's inequality with \eqref{eq:torus-L2}, \eqref{eq:joint-embedding} and \eqref{eq:torus-Hs-cor}, we obtain 
\begin{align*}
\norm{\boldb(t)}{H^{s}}^{2}+\int^{t}_{0}\norm{\boldu(\tau)}{\dot{H}^{s+\alpha}}^{2}\myd{\tau}\leq \norm{\boldb_{0}}{H^{s}}^{2}\exp{\left(C\int^{t}_{0}\norm{\nabla \boldu(\tau)}{L^{\infty}}+{\norm{\nabla\boldb(\tau)}{L^{q}}^{2}}\myd{\tau}\right)},
\end{align*}
{where $2<q\leq\infty$ is defined in \eqref{eq:p-q}. {We note that the quantity $\norm{\nabla\boldb}{L^{\infty}}$ in the bound in {\cite[Theorem 2.2]{BFV22}} is replaced by $\norm{\nabla\boldb}{L^{q}}$.}

To show $\boldu\in C([0,T);\thSob{s+2\alpha-\frac{d}{2}})$, we use the argument used to derive \eqref{eq:u-estimate-b} under the condition (\ref{LWP Condition 2}):
\[
\begin{aligned}
\norm{\Lambda^{s+2\alpha-\frac{d}{2}}\boldu(t)}{L^{2}}&=\norm{\Lambda^{s-\frac{d}{2}}\mathbb{P}\Div(\boldb\otimes\boldb)(t)}{L^{2}}\\
&\apprle \norm{\Lambda^{s+1-\frac{d}{2}}\left(\boldb(t)\otimes\boldb(t)\right)}{L^{2}} \\
&\apprle \norm{\Lambda^{s+1-\frac{d}{2}}\boldb(t)}{\Leb{d}}\norm{\boldb(t)}{\Leb{r}} \apprle \norm{\boldb(t)}{H^{s}}^{2},
\end{aligned}
\]
where $r=6$ if $d=3$ and $r=\infty$ if $d=2$. This completes the proof of Theorem \ref{thm:LWP}.

%%%%%%%%%%%%%%%%%%%%%%%%
\subsection{Proof of Theorem \ref{thm:A}}
%%%%%%%%%%%%%%%%%%%%%%
Suppose  $\alpha>\frac{d}{2}$ and $s\geq 1$. By Proposition \ref{prop:BKM-estimate} with $\alpha>\frac{d}{2}$, we have
\begin{equation}\label{eq:BKM-estimate-2}
\norm{\nabla \boldu}{\Leb{\infty}}\apprle 1+\norm{\boldu}{\thSob{\frac{d}{2}+1}}\log (e+\norm{\boldu}{\thSob{1+\alpha}}).
\end{equation}
By \eqref{eq:torus-L2}, \eqref{eq:torus-H1}, \eqref{eq:BKM-estimate-2}, and \eqref{eq:u-estimate-b} with $s=1$, we obtain 
\begin{equation}\label{eq:torus-H1-log}
\begin{aligned}
\frac{1}{2}\frac{d}{dt}\norm{\boldb}{H^{1}}^{2}+\norm{\boldu}{\dot{H}^{1+\alpha}}^{2} &\leq C\norm{\boldb}{\tSob{1}}^2+ C\norm{\boldu}{\dot{H}^{\frac{d}{2}+1}}\log (e+\norm{\boldu}{\thSob{1+\alpha}})\norm{\boldb}{H^{1}}^{2} \\
&\leq C\norm{\boldb}{\tSob{1}}^2+C\norm{\boldu}{\dot{H}^{\frac{d}{2}+1}}\norm{\boldb}{H^{1}}^{2}\log \left(e+\norm{\boldb}{\tSob{1}}^2\right).
\end{aligned}
\end{equation}
Define 
\[
y(t)=e+\norm{\boldb(t)}{\tSob{1}}^2+2\int_0^t \norm{\boldu(\tau)}{\thSob{1+\alpha}}^2 \myd{\tau}.
\]
Then by \eqref{eq:torus-H1-log}, we have
\begin{align*}
\frac{d}{dt}y(t)&\leq C\norm{\boldu(t)}{\dot{H}^{\frac{d}{2}+1}}y(t)\log y(t)+Cy(t).
\end{align*}
Hence by dividing $y(t)$ and using Gr\"onwall's inequality, we get
\begin{equation}\label{eq:gronwall-final}
 \log y(t) \leq [\log y(0)+Ct]\exp\left( C\int_0^t \norm{\boldu(\tau)}{\dot{H}^{\frac{d}{2}+1}}d\tau\right).
\end{equation}
{Therefore by \eqref{eq:gronwall-final}, \eqref{eq:BKM}, and \eqref{eq:torus-exp} with $\norm{\nabla\boldu}{L^{\infty}}\apprle\norm{\boldu}{\dot{H}^{1+\alpha}}$,} we get the desired result Theorem \ref{thm:A} (i). 

Theorem \ref{thm:A} (ii) is an immediate consequence of Theorem \ref{thm:A} (i) because \eqref{eq:torus-energy-ineq} {with $\alpha\geq\frac{d}{2}+1$} implies that 
\[
\int_0^T \norm{\boldu(t)}{\thSob{\frac{d}{2}+1}}dt<\infty
\]
for all $T>0$. This completes the proof of Theorem \ref{thm:A}.

%%%%%%%%%%%%%%%%%%%%%%%%%%%%%%%%%%
\section{Resistive case: Proof of Theorem \ref{thm:B}}\label{sec:5}
%%%%%%%%%%%%%%%%%%%%%%%%%%%%%%%%%%
We recall \eqref{eq:MRE} with $\eta=1$ {and $\nu=1$}:
\begin{equation}\label{eq:MRE 5}
\left\{
\begin{alignedat}{2}
\Lambda^{2\alpha} \boldu +\nabla \pbar &=(\boldb\cdot \nabla)\boldb&&\quad \text{in } \mathbb{R}^d\times (0,T),\\
\partial_t \boldb+\Lambda^{2\beta} \boldb +(\boldu\cdot \nabla)\boldb &=(\boldb\cdot \nabla)\boldu &&\quad \text{in }\mathbb{R}^d \times(0,T),\\
\Div \boldu =\Div \boldb &=0 &&\quad \text{in } \mathbb{R}^d\times (0,T),\\
\boldb(\cdot,0)&=\boldb_0 &&\quad \text{on } \mathbb{R}^d.
\end{alignedat}\right.
\end{equation}

We define 
\[    
\frac{1}{q}:=\frac{1}{p}-\frac{2\beta-1}{3d}=\frac{3\alpha+\beta-2}{3d}, \quad \frac{1}{r}:=\frac{2}{q}-\frac{2\alpha-1}{d}.
\]
By the condition of $(\alpha,\beta)$ in Theorem \ref{thm:B}, it is easy to check that 
\begin{equation}\label{eq:p-q-r}
   1<\max\{p,2\}<q<\infty\quad \text{and}\quad 0<\frac{1}{r}=\frac{2\beta-1}{3d}<1.
\end{equation} 
Moreover, we have
\[  \frac{1}{p}=\frac{1}{r}+\frac{1}{q}.\]
For these $p$ and $q$, we introduce the following norms: 
\[
\begin{split}
\norm{\boldb}{\mathcal{N}_T}&=\norm{\boldb}{\Leb{\infty}(0,T;\Leb{p})},\quad \norm{\boldb}{\mathcal{F}_T}=\sup_{0<s\leq T} s^{\frac{d}{2\beta}\left(\frac{1}{p}-\frac{1}{q}\right)} \norm{\boldb(s)}{\Leb{q}},\\
\mathcal{I}_T&=\sup_{0<t\leq T} t^{\frac{d}{2\beta}\left(\frac{1}{p}-\frac{1}{q}\right)}\norm{G_{\beta}(t)\boldb_0}{\Leb{q}}.
\end{split} 
\]

Since $0<2\alpha-1<d$ and $2<q<\infty$, it follows from Hardy-Littlewood-Sobolev inequality (Proposition  \ref{thm:HLS}) that
\begin{equation}\label{eq:u-r-estimate}
\begin{aligned}
   \norm{\boldu}{\Leb{r}}&\approx \norm{\Lambda^{1-2\alpha} \Proj(\boldb\otimes\boldb)}{\Leb{r}}\apprle \norm{\Proj (\boldb\otimes\boldb)}{\Leb{\frac{q}{2}}}\\
   &\apprle \norm{\boldb \otimes\boldb}{\Leb{\frac{q}{2}}}\apprle \norm{\boldb}{\Leb{q}}^2.
\end{aligned}
\end{equation}
Hence by H\"older's inequality and \eqref{eq:u-r-estimate}, we {have}
\begin{equation}\label{eq:holder}
     \norm{\boldu\otimes\boldb}{\Leb{p}}\apprle \norm{\boldu}{\Leb{r}}\norm{\boldb}{\Leb{q}}\apprle \norm{\boldb}{\Leb{q}}^3.
\end{equation}

We now define the operator $\mathcal{T}$: 
\begin{equation}\label{eq:mild-solution-definition}
    \mathcal{T}(\boldb,\boldb_0)(t)= G_{\beta}(t)\boldb_0+\mathcal{N}(\boldb)(t),\quad (\boldb,\boldb_0)\in \mathcal{F}_{T_*} \times \Leb{p}_\sigma,
\end{equation}
where
\[ 
\mathcal{N}(\boldb)(t)=-\int_0^t G_{\beta}(t-s) \Div(\boldu\otimes\boldb-\boldb\otimes \boldu)(s)\myd{s}, \quad \boldu=\Lambda^{-2\alpha} \Proj\Div(\boldb\otimes\boldb).
\]
Let 
\[  
\mathcal{B}_{\varepsilon,T_*}=\{ \boldb \in \mathcal{F}_{T_*} : \norm{\boldb}{\mathcal{F}_{T_*}}\leq \varepsilon,\quad \norm{\boldb}{\mathcal{N}_{T_*}} \leq \norm{\boldb_0}{\Leb{p}}+\varepsilon\}.
 \]

\vspace{1ex}

For fixed $\boldb_0 \in \Leb{p}_\sigma$, we first show that $\mathcal{T}(\cdot,\boldb_0)$ maps from $\mathcal{B}_{\varepsilon_2,T_2}$ to itself for small $T_2$ and $\varepsilon_2$.  By Proposition \ref{prop:linear-estimate} and \eqref{eq:holder}, we first bound $\mathcal{T}(\boldb,\boldb_0)(t)$ in $L^{p}$: 
\begin{align*}
\norm{\mathcal{T}(\boldb,\boldb_0)(t)}{\Leb{p}}&\leq \norm{G_{\beta}(t) \boldb_0}{\Leb{p}}+2C\int_0^t (t-s)^{-\frac{1}{2\beta}} \norm{(\boldu\otimes\boldb)(s)}{\Leb{p}}\myd{s}\\
&\leq \norm{G_{\beta}(t) \boldb_0}{\Leb{p}}+2C\int_0^t (t-s)^{-\frac{1}{2\beta}} \norm{\boldb(s)}{\Leb{q}}^3\myd{s}.
\end{align*}
Let
\[
\sigma=\frac{d}{2\beta}\left(\frac{1}{p}-\frac{1}{q}\right) \quad \text{and so} \quad 3\sigma+\frac{1}{2\beta}=1.
\]
Then by \eqref{eq:beta-function-estimate},  we have
\begin{align*}
\norm{\mathcal{T}(\boldb,\boldb_0)(t)}{\Leb{p}}&\leq \norm{G_{\beta}(t) \boldb_0}{\Leb{p}}+2C\int_0^t (t-s)^{-\frac{1}{2\beta}} s^{-3\sigma} s^{3\sigma}\norm{\boldb(s)}{\Leb{q}}^3\myd{s}\\
&\leq \norm{G_{\beta}(t) \boldb_0}{\Leb{p}}+ 2C\sup_{0<s<t} (s^\sigma\norm{\boldb(s)}{\Leb{q}})^3 \int_0^t (t-s)^{-\frac{1}{2\beta}} s^{-3\sigma} \myd{s}\\
&\leq \norm{G_{\beta}(t) \boldb_0}{\Leb{p}}+ 2C\sup_{0<s<t} (s^\sigma\norm{\boldb(s)}{\Leb{q}})^3.
\end{align*} 
Hence by Proposition \ref{prop:linear-estimate}, we have
\begin{equation}\label{eq:Nt-mapping}
\norm{\mathcal{T}(\boldb,\boldb_0)}{\mathcal{N}_T} \leq C_1 \norm{\boldb_0}{\Leb{p}}+ C_2 \norm{\boldb}{\mathcal{F}_T}^3
\end{equation}
for some constants $C_1,C_2>0$. 
     
Next we estimate \eqref{eq:mild-solution-definition} in $\Leb{q}$-norm. By $2\beta>1$ and \eqref{eq:p-q-r}, we have
\eqn \label{eq:beta-sigma-sum}
0<\frac{1}{2\beta}+\sigma<\frac{1}{2\beta}+3\sigma=1, \quad 1-\left(\frac{1}{2\beta}+\sigma\right)-3\sigma=-\sigma.
\een 
It follows from Proposition \ref{prop:linear-estimate},  \eqref{eq:p-q-r}, \eqref{eq:holder}, and \eqref{eq:beta-function-estimate} that 
\begin{align}
\norm{\mathcal{T}(\boldb,\boldb_0)(t)}{\Leb{q}}&\leq \norm{G_{\beta}(t) \boldb_0}{\Leb{q}}+2\int_0^t \norm{G_{\beta}(t-s)   \Div(\boldu\otimes\boldb)(s)}{\Leb{q}}\myd{s}\label{eq:Ft-estimate}\\
&\leq \norm{G_{\beta}(t) \boldb_0}{\Leb{q}}+C\int_0^t (t-s)^{-\frac{1}{2\beta}-\sigma}\norm{(\boldu\otimes\boldb)(s)}{\Leb{p}}\myd{s}\nonumber\\
&\leq \norm{G_{\beta}(t) \boldb_0}{\Leb{q}}+ C\sup_{0<s<t} (s^\sigma\norm{\boldb(s)}{\Leb{q}})^3 \int_0^t  (t-s)^{-\frac{1}{2\beta}-\sigma}s^{-3\sigma} \myd{s}\nonumber\\
&\leq \norm{G_{\beta}(t) \boldb_0}{\Leb{q}}+Ct^{-\sigma} \sup_{0<s<t} (s^\sigma\norm{\boldb(s)}{\Leb{q}})^3,\nonumber
\end{align}
where we use \eqref{eq:beta-sigma-sum} to the last equality. Multiplying $t^\sigma$ to this inequality and taking supremum over $0<t<T$, we deduce that 
\begin{equation}\label{eq:Ft-mapping}
\norm{\mathcal{T}(\boldb,\boldb_0)}{\mathcal{F}_T}\leq \mathcal{I}_T +C_3 \norm{\boldb}{\mathcal{F}_T}^3.
\end{equation}
Choose $0<\varepsilon_1<\min\{1,\frac{1}{4C_3}\}$. By Proposition \ref{prop:time-control}, there exists $T_1>0$ such that $\mathcal{I}_{T_1}\leq \frac{\varepsilon}{2}$. Hence 
\[
   \norm{\mathcal{T}(\boldb,\boldb_0)}{\mathcal{F}_{T_1}}\leq \frac{\varepsilon_1}{2}+\frac{\varepsilon_1^2}{4}\leq \varepsilon_1.
\]
Moreover, we have
\[
   \norm{\mathcal{T}(\boldb,\boldb_0)(t)}{\Leb{p}}\leq \norm{G_{\beta}(t)\boldb_0}{\Leb{p}}+C_4\varepsilon^3_1,
\]
for all $t \in [0,T_1]$. Choose $0<\varepsilon_2\leq \min\{\varepsilon_1, \frac{1}{2C_4}\}$. Since $G_{\beta}(t):\Leb{p}\rightarrow \Leb{p}$ is continuous in $t$, we can choose $T_2<T_1$ sufficiently small so that 
\[
\sup_{t\in [0,T_2]} \norm{G_{\beta}(t)\boldb_0}{\Leb{p}}\leq \norm{\boldb_0}{\Leb{p}}+\frac{\varepsilon_2}{2}.
\]
Therefore, $\mathcal{T}(\cdot,\boldb_0)$ maps from $\mathcal{B}_{\varepsilon_2,T_2}$ to itself for small $T_2$ and $\varepsilon_2$. 

Next, we show that $\mathcal{T}$ is a uniform contraction on $\mathcal{B}_{\varepsilon,{T_*}}$ for sufficiently small $0<\varepsilon<\varepsilon_2$ and $0<{T_*}<T_2$. For $\boldb_i \in \mathcal{F}_T$ and $\boldu_i = \Lambda^{-2\alpha} \Proj\Div(\boldb_i\otimes\boldb_i)$, $i=1,2$, if we write $\widetilde{\boldb}=\boldb_1-\boldb_2$ and $\widetilde{\boldu}=\boldu_1-\boldu_2$, then 
\[
\widetilde{\boldu}=\Lambda^{-2\alpha}\Proj\Div(\widetilde{\boldb}\otimes\boldb_1)+\Lambda^{-2\alpha}\Proj\Div(\boldb_2\otimes\widetilde{\boldb})
\]
and
\begin{align*}
&\relphantom{=}\mathcal{T}(\boldb_1,\boldb_0)(t)-\mathcal{T}(\boldb_2,\boldb_0)(t)=\mathcal{N}(\boldb_1)(t)-\mathcal{N}(\boldb_2)(t)\\
&=\int_0^t G_{\beta}(t-s) \Div(\tilde{\boldb}\otimes\boldu_1+\boldb_2\otimes\tilde{\boldu}-\boldu_1\otimes\tilde{\boldb}-\tilde{\boldu}\otimes\boldb_2)(s)\myd{s}.
\end{align*}
Following the same arguments in \eqref{eq:u-r-estimate} and \eqref{eq:Ft-mapping}, we get  
\begin{equation*}
\norm{\widetilde{\boldu}}{\Leb{r}}\apprle \left(\norm{\boldb_1}{\Leb{q}}+\norm{\boldb_2}{\Leb{q}}\right)\norm{\widetilde{\boldb}}{\Leb{q}}.
\end{equation*}
and
\begin{align*}
&\norm{\mathcal{N}(\boldb_1)(t)-\mathcal{N}(\boldb_2)(t)}{\Leb{q}}\\
&\leq C\int_0^t (t-s)^{-\frac{1}{2\beta}-\sigma}\left(\norm{\widetilde{\boldb}\otimes {\boldu_1}}{\Leb{p}}+\norm{{\boldb}_2\otimes \widetilde{\boldu}}{\Leb{p}}+\norm{\boldu_1\otimes\widetilde{\boldb}}{\Leb{p}}+\norm{\widetilde{\boldu}\otimes\boldb_2}{\Leb{p}} \right)\myd{s}\\
&\leq C t^{-\sigma}\left(\norm{\widetilde{\boldb}}{\mathcal{F}_T}\norm{{\boldb}_1}{\mathcal{F}_T}^2+\norm{{\boldb}_2}{\mathcal{F}_T}\left(\norm{\boldb_1}{\mathcal{F}_T}+\norm{\boldb_2}{\mathcal{F}_T}\right)\norm{\tilde{\boldb}}{\mathcal{F}_T}\right).
\end{align*}
This implies that there exists a constant $C=C(d,\alpha,\beta)>0$ uniformly in $\boldb_0$ such that 
\begin{equation}\label{eq:uniform-contraction}
   \norm{\mathcal{T}(\boldb_1,\boldb_0)-\mathcal{T}(\boldb_2,\boldb_0)}{{\mathcal{F}_{T_*}}}\leq C \varepsilon^2\norm{\boldb_1-\boldb_2}{{\mathcal{F}_{T_*}}}
\end{equation}
for all $\boldb_i \in B_{\varepsilon,T_*}$. By choosing $\varepsilon>0$  so that $C\varepsilon^2<\frac{1}{2}$, we see that $\mathcal{T}$ is a contractive mapping from $\mathcal{B}_{\varepsilon,{T_*}}$ onto itself. By the contraction principle, there exists $g:\Leb{p}\rightarrow \mathcal{B}_{\varepsilon,{T_*}}\subset \mathcal{F}_{{T_*}}$ such that  
\[  
g(\boldb_0)=\mathcal{T}(g(\boldb_0),\boldb_0).
\]

 In order to show the continuous dependance of $g(\boldb_0)$, which also gives the uniqueness of solutions, we take two initial data $\boldb_0^i \in \Leb{p}_\sigma(\mathbb{R}^d)$, $i=1,2$. Then 
\[
g(\boldb_0^1)-g(\boldb_0^2)=G_\beta(t)\left(\boldb_0^1-\boldb_0^2\right)+[\mathcal{N}(g(\boldb_0^1))-\mathcal{N}(g(\boldb_0^2))](t).
\]
 By \eqref{eq:uniform-contraction} and Proposition \ref{prop:linear-estimate}, we obtain
\[
\norm{g(\boldb_0^1)-g(\boldb_0^2)}{\mathcal{F}_{{T_*}}}\leq C_1\norm{\boldb_0^1-\boldb_0^2}{\Leb{p}}+C_2\varepsilon^2\norm{g(\boldb_0^1)-g(\boldb_0^2)}{\mathcal{F}_{{T_*}}},
\]
where $C_2$ is the same constant $C$ in \eqref{eq:uniform-contraction}. This implies
\begin{equation}\label{eq:mild-solution-difference}
\norm{g(\boldb_0^1)-g(\boldb_0^2)}{\mathcal{F}_{{T_*}}}\leq C(d,\alpha,\beta,p) \norm{\boldb_0^1-\boldb_0^2}{\Leb{p}}.
\end{equation}
Hence by following the proof of \eqref{eq:Nt-mapping} and using \eqref{eq:mild-solution-difference}, we obtain
\begin{equation}\label{eq:mild-solution-difference-2}
\norm{g(\boldb_0^1)-g(\boldb_0^2)}{\mathcal{N}_{{T_*}}}\leq C\norm{\boldb_0^1-\boldb_0^2}{\Leb{p}}+ C\varepsilon^2\norm{g(\boldb_0^1)-g(\boldb_0^2)}{\mathcal{F}_{T_*}}\apprle \norm{\boldb_0^1-\boldb_0^2}{\Leb{p}}.
\end{equation}

By a standard argument, one can show that $\boldb \in C([0,T_*];\Leb{p})$ and mild solutions are unique in the class $C([0,T];\Leb{p})\cap \mathcal{F}_T$ (see e.g. \cite[Proof of Theorem 5.12]{BV22}).  Also, following the argument in the above, one can see that there exists $\varepsilon_0>0$ such that if $\norm{\boldb_0}{\Leb{p}}<\varepsilon_0$, then there exists a unique global mild solution $\boldb$ to {\eqref{eq:MRE 5}}. 

\vspace{1ex}

To complete the proof of Theorem \ref{thm:B}, it remains  to show the asymptotic behavior of solutions. We  observe below that if $\boldb_0$ also belongs to $\Leb{s}$ for some $1\leq s<p$, then we can obtain temporal decay rate for the norm $\norm{\boldb(t)}{\Leb{p}}$.

\begin{prop}\label{prop:convergence}
Let $\varepsilon_0>0$ be a number so that the mild solution $\boldb$ constructed in Theorem \ref{thm:B} globally exists for  $\boldb_0$ satisfying $\norm{\boldb_0}{\Leb{p}}<\varepsilon_0$. Then there exist $s \in [1,p)$ and $0<\varepsilon_1\leq \varepsilon_0$ such that if $\boldb_0\in \Lebdiv{s}(\mathbb{R}^d)\cap \Lebdiv{p}(\mathbb{R}^d)$ satisfies $\norm{\boldb_0}{\Leb{p}}<\varepsilon_1$, then 
\[       
t^{\theta} \boldb(t)\in BC([0,\infty);\Leb{p}), \quad \theta=\frac{d}{2\beta}\left(\frac{1}{s}-\frac{1}{p}\right).
\]
\end{prop}

\begin{proof}
We first notice  that 
\[     
\frac{2\beta}{d}+2\left(\frac{1}{q}-\frac{1}{p}\right)=\frac{2\beta}{d}+\frac{2(1-2\beta)}{3d}=\frac{2\beta+2}{3d}>0.
\]
Hence we can choose $1\leq s<p<\infty$ so that 
\begin{equation}\label{eq:s-p}
    \frac{1}{p}<\frac{1}{s}<\frac{2\beta}{d}+\frac{2}{q}-\frac{1}{p}.
\end{equation}
Let 
\[    
\theta = \frac{d}{2\beta}\left(\frac{1}{s}-\frac{1}{p}\right)\quad\text{and}\quad \sigma=\frac{d}{2\beta}\left(\frac{1}{p}-\frac{1}{q}\right).
\]
By \eqref{eq:s-p}, we have
\[
  \theta>0,\quad 0< 2\sigma+\theta <1.
\] 
We now introduce 
\[
\norm{\boldb}{Z_T}=\sup_{0<\tau<T} \left(\tau^\theta \norm{\boldb(\tau)}{\Leb{p}}\right).
\]
Since $1/(2\beta)+3\sigma=1$ and $0<2\sigma+\theta<1$, 
\begin{align*}
\norm{\boldb(t)}{\Leb{p}}&\apprle t^{-\theta}\norm{\boldb_0}{\Leb{s}}+\int_0^t (t-\tau)^{-\frac{1}{2\beta}-\frac{d}{2\beta r}} \norm{\boldu(\tau)}{\Leb{r}}\norm{\boldb(\tau)}{\Leb{p}}\myd{\tau}\\
&\apprle t^{-\theta}\norm{\boldb_0}{\Leb{s}}+\int_0^t (t-\tau)^{-\frac{1}{2\beta}-\sigma}\tau^{-2\sigma} \norm{\boldb}{\mathcal{F}_\infty}^2 \tau^{-\theta} \tau^\theta \norm{\boldb(\tau)}{\Leb{p}}\myd{\tau}\\
&\apprle t^{-\theta} \norm{\boldb_0}{\Leb{s}}+t^{-\theta} \norm{\boldb}{\mathcal{F}_\infty}^2 {\sup_{0<\tau<t}\left(\tau^{\theta}\norm{\boldb(\tau)}{L^{p}}\right)}.
\end{align*}
Multiplying $t^\theta$ to this inequality, we get
\[      
\norm{\boldb}{Z_T}\apprle \norm{\boldb_0}{\Leb{s}}+\norm{\boldb}{\mathcal{F}_\infty}^2\norm{\boldb}{Z_T},
\]
where the implicit constant does not depend on $T$. Then by \eqref{eq:Ft-mapping} and choosing $\varepsilon_1$ sufficiently small, we obtain  the desired result.
\end{proof}

Now we are ready to prove the asymptotic behavior of mild solutions. Let $\varepsilon>0$ be given and let $\varepsilon_0>0$ be a small constant  so that for $\norm{\boldb_0}{\Leb{p}}<\varepsilon_0$, there exists a global solution associated with $\boldb_0$.  For such $\boldb_0^1$ and $\boldb_0^2$, it follows from {\eqref{eq:mild-solution-difference-2}} that if $\boldb_i$ is a unique mild solution with the initial data $\boldb_0^i$, $i=1,2$, respectively, then there exists a constant $C_0$ depends on $d$, $\alpha$, $\beta$, and $\varepsilon_0$
\begin{equation}\label{eq:boldb-mapping}
    \sup_{t>0} \norm{\boldb_1(t)-\boldb_2(t)}{\Leb{p}}\leq C_0\norm{\boldb_0^1-\boldb_0^2}{\Leb{p}}.
\end{equation}
For $\boldb_0 \in \Lebdiv{p}$ satisfying $\norm{\boldb_0}{\Leb{p}}\leq \frac{\varepsilon}{2}$, where $0<{\varepsilon}<\varepsilon_0$ is given in Proposition \ref{prop:convergence}, we choose $\widetilde{\boldb}_0 \in C_{c,\sigma}^\infty$ so that 
\begin{equation}\label{eq:b-b0-tilde}
\norm{\widetilde{\boldb}_0-\boldb_0}{\Leb{p}}\leq \min\left\{ \varepsilon/C_0,{\varepsilon}/{2}\right\},
\end{equation}
 where $C_0$ is the same constant in \eqref{eq:boldb-mapping}. Then $\norm{\widetilde{\boldb}_0}{\Leb{p}}\leq {\varepsilon}$. Let $\widetilde{\boldb}$ be the global mild solution with the initial data $\widetilde{\boldb}_0$. By choosing $\varepsilon>0$ sufficiently small, \eqref{eq:boldb-mapping}, and \eqref{eq:b-b0-tilde}, we have
 \[
     \norm{\boldb(t)-\widetilde{\boldb}(t)}{\Leb{p}}\leq \varepsilon
\]
for all $t>0$. On the other hand, it follows from Proposition \ref{prop:convergence} that 
\[
\lim_{t\rightarrow\infty} \norm{\widetilde{\boldb}(t)}{\Leb{p}}=0
\]
which implies that 
\[
\norm{\boldb(t)}{\Leb{p}}\leq 2\varepsilon
 \]
for sufficiently large $t$. Hence 
\[
\limsup_{t\rightarrow\infty} \norm{\boldb(t)}{\Leb{p}}\leq 2\varepsilon.
\]
Since $\varepsilon>0$ is arbitrary and  $\norm{\boldb(t)}{\Leb{p}}\ge 0$, we arrive at the desired result and this completes the proof of Theorem \ref{thm:B}.

%%%%%%%%%%%%%%%
\section*{Acknowledgments}
%%%%%%%%%%%%%%%

H. Bae and J. Shin were supported by the National Research Foundation of Korea(NRF) grant funded by the Korea government(MSIT) (grant No. 2022R1A4A1032094). H. Kwon was partially supported by the NSF under agreement DMS-2055244 and the international travel fund award by Brown University.

\bibliographystyle{amsplain}

\begin{thebibliography}{99}


\bibitem{Arn74}
V.~I. Arnold. 
	\newblock The asymptotic Hopf invariant and its applications.
	\newblock \emph{Selecta Math. Soviet.} \textbf{5} (1986), no.~4, 327--345, Selected
  translations. \MR{891881}
  
  

\bibitem{BBT12}
H. Bae, A. Biswas, E. Tadmor.
	\newblock Analyticity and decay estimates of the Navier-Stokes equations in critical Besov spaces.
	\newblock \emph{Arch. Ration. Mech. Anal.} {\bf 205} (2012), no. 3, 963--991. \MR{2960037}



\bibitem{BKM84}
J.~T. Beale, T. Kato, A. Majda.
	\newblock Remarks on the breakdown of smooth solutions for the 3-D Euler equations.
	\newblock \emph{Comm. Math. Phys.} {\bf 94} (1984), no. 1, 61--66. \MR{763762}



\bibitem{BV22}
J. Bedrossian, V. Vicol.
	\newblock The mathematical analysis of the incompressible Euler and Navier-Stokes equations - an introduction.
	\newblock Graduate Studies in Mathematics, 225. \emph{American Mathematical Society, Providence, RI}, [2022], xiii+218 pp. \MR{4475666}



\bibitem{BFV22}
R. Beekie, S. Friedlander, V. Vicol.
     \newblock On Moffatt's magnetic relaxation equations.
     \newblock \emph{Comm. Math. Phys.} {\bf 390} (2022) no. 3, 1311--1339.



\bibitem{BO13}
\'A. B\'enyi, T.  Oh.
	\newblock The Sobolev inequality on the   torus revisited.
	\newblock \emph{Publ. Math. Debrecen} {\bf83} (2013), no. 3, 359--374.
  \MR{3119672}



\bibitem{Bre14}
Y. Brenier.
	\newblock Topology-preserving diffusion of divergence-free vector fields and magnetic relaxation.
	\newblock \emph{Comm. Math. Phys.} {\bf 330} (2014), no. 2, 757--770.



\bibitem{BG80}
	H. Br\'{e}zis, T. Gallouet.
	\newblock Nonlinear Schr\"{o}dinger evolution equations.
	\newblock \emph{Nonlinear Anal}. {\bf  4} (1980), no. 4, 677–681. \MR{582536}



\bibitem{BW80}
	\newblock H. Br\'{e}zis, S. Wainger.
	\newblock A note on limiting cases of Sobolev embeddings and convolution inequalities.
	\newblock \emph{Comm. Partial Differential Equations} {\bf 5} (1980), no. 7, 773--789. \MR{579997}



\bibitem{BL96}
	O. P. Bruno, Peter, Laurence.
	\newblock Existence of three-dimensional toroidal MHD equilibria with nonconstant pressure.
	\newblock \emph{Comm. Pure Appl. Math.} {\bf49} (1996), no. 7, 717--764. \MR{1387191}



\bibitem{CK19}
D. Cardona, V.  Kumar.
	\newblock $L^p$-boundedness and $L^p$-nuclearity of multilinear pseudo-differential operators on ${\Bbb  Z}^n$ and the torus ${\Bbb T}^n$.
	\newblock \emph{J. Fourier Anal. Appl.} {\bf25} (2019), no. 6, 2973--3017. \MR{4029168}



\bibitem{CF08}
J. A. Carrillo, L.C. F. Ferreira.
	\newblock The asymptotic behaviour  of subcritical dissipative quasi-geostrophic equations.
	\newblock \emph{Nonlinearity} \textbf{21} (2008), no.~5, 1001--1018. \MR{2412324}



\bibitem{C02}
D. Chae.
	\newblock On the well-posedness of the Euler equations in the Triebel-Lizorkin spaces.
	\newblock \emph{Comm. Pure Appl. Math.} \textbf{55} (2002),  no.~5, 654--678. \MR{1880646}

\bibitem{CG19}
M. Cirant, A. Goffi.
	\newblock On the existence and uniqueness of solutions to time-dependent fractional MFG.
	\newblock \emph{SIAM J. Math. Anal.} \textbf{51} (2019), no.~2, 913--954. \MR{3934106}



\bibitem{CP22}
P. Constantin, F. Pasqualotto.
	\newblock Magnetic relaxation of a Voigt-MHD system.
	\newblock arXiv:2208.11109.

\bibitem{D01}
P. A. Davidson.
\newblock An Introduction to Magnetohydrodynamics,
\newblock \emph{Cambridge Texts Appl. Math., Cambridge University Press, Cambridge,} 2001.

\bibitem{ELP21}
A. Enciso, A. Luque, D. Peralta-Salas.
	\newblock MHD equilibria with nonconstant pressure in nondegenerate toroidal domains.
	\newblock To appear in \emph{J. Eur. Math. Soc. (JEMS)}.



\bibitem{Evans}
L. Evans.
\newblock Partial differential equations. Graduate Studies in Mathematics, 19.
\newblock \emph{American Mathematical Society, Providence, RI,} 1998. xviii+662 pp. \MR{1625845}




\bibitem{FMRR14}
C.L. Fefferman, D.S. McCormick, J.C. Robinson, J. L. Rodrigo.
	\newblock Higher order commutator estimates and local existence for the  non-resistive MHD equations and related models.
	\newblock  \emph{J. Funct. Anal.}  \textbf{267} (2014), no.4, 1035--1056. \MR{3217057}



\bibitem{Feng23}
W. Feng, F. Hafeez, J.  Wu, D. Regmi.
	\newblock Stability and   exponential decay for magnetohydrodynamic equations.
	\newblock \emph{Proc. Roy. Soc. Edinburgh Sect. A} {\bf153} (2023), no. 3, 853--880. \MR{4595825}



\bibitem{G67}
H. Grad.
	\newblock Toroidal Containment of a Plasma.
	\newblock \emph{Phys. Fluids.}  \textbf{10} (1967), 137--154.



\bibitem{JT21}
Y. Ji, W. Tan.
	\newblock Global well-posedness of a 3D Stokes-Magneto equations with fractional magnetic diffusion.
  	\newblock \emph{Discrete Contin. Dyn. Syst. Ser. B} {\bf26} (2021), no. 6, 3271--3278. \MR{4235654}



\bibitem{JT21-1}
Y. Ji, W. Tan.
	\newblock Large time behavior of solutions to a Stokes-Magneto equations in three dimensions.
  	\newblock \emph{J. Evol. Equ.} {\bf 21} (2021), no. 2, 2449--2470. \MR{4278435}



\bibitem{Ju04}
N. Ju.
	\newblock Existence and uniqueness of the solution to the dissipative 2D quasi-geostrophic equations in the Sobolev space.
  \newblock \emph{Comm. Math. Phys.} {\bf 251} (2004), no. 2, 365--376. \MR{2100059}



\bibitem{KP88}
T. Kato, G. Ponce.
	\newblock Commutator estimates and the Euler and Navier-Stokes equations.
  	\newblock \emph{Comm. Pure Appl. Math.} {\bf 41} (1988), no. 7, 891--907. \MR{951744}



\bibitem{KPV91}
C. E. Kenig, G. Ponce, L. Vega.
	\newblock Well-posedness of the  initial value problem for the Korteweg-de Vries equation.
	\newblock \emph{J. Amer. Math. Soc.} {\bf 4} (1991), no. 2, 323--347.  \MR{1086966}




\bibitem{KK23}
H. Kim, H. Kwon.
     \newblock Global existence and uniqueness of weak solutions of a Stokes-Magneto system with fractional diffusions.
     \newblock \emph{J. Differential Equations} {\bf 374} (2023) 497--547. \MR{4626420}





\bibitem{KOT02}
H. Kozono, T.  Ogawa, Y. Taniuchi.
	\newblock The critical  Sobolev inequalities in Besov spaces and regularity criterion to some
  semi-linear evolution equations.
 	 \newblock \emph{Math. Z.} {\bf 242} (2002), no. 2, 251--278. \MR{1980623}

\bibitem{KT00}
H. Kozono. Y. Taniuchi.
	\newblock Limiting case of the Sobolev  inequality in BMO, with application to the Euler equations.
  	\newblock \emph{Comm. Math. Phys.} {\bf 214} (2000), no. 1, 191--200. \MR{1794270}


\bibitem{L19}
D. Li.
	\newblock On Kato-Ponce and fractional Leibniz.
	\newblock \emph{Rev. Mat. Iberoam.} {\bf 35} (2019), no. 1, 23--100. \MR{3914540}



\bibitem{Lieb01}
E. H. Lieb, M. Michael Loss.
	\newblock Analysis.
	\newblock Second edition. Graduate Studies in Mathematics, 14. \emph{American Mathematical Society, Providence, RI,} 2001. xxii+346 pp.  \MR{1415616}


\bibitem{MRR14}
D. S. McCormick, J. C. Robinson, J. L. Rodrigo.
	\newblock Existence and   uniqueness for a coupled parabolic-elliptic model with applications to
  magnetic relaxation.
	  \newblock \emph{Arch. Ration. Mech. Anal.} {\bf 214} (2014), no. 2, 503--523. \MR{3255698}


\bibitem{MYZ08}
C. Miao, B. Yuan, B. Zhang.
	\newblock Well-posedness of the Cauchy problem for the fractional power dissipative equations.
	\newblock \emph{Nonlinear Anal.} {\bf 68} (2008), no. 3, 461--484.  \MR{2372358}



\bibitem{M85}
H. K. Moffatt.
	\newblock Magnetostatic equilibria and analogous Euler flows of arbitrarily complex topology. {I}. Fundamentals.
 	 \newblock \emph{J. Fluid Mech.} {\bf159} (1985), 359--378. \MR{819398}



\bibitem{M21}
H. K. Moffatt. 
	\newblock Some topological aspects of fluid dynamics.
	\newblock \emph{J. Fluid Mech.} {\bf 914} (2021), Paper No. P1, 56 pp. \MR{4232242}



\bibitem{O03}
T.  Ogawa.
	\newblock Sharp Sobolev inequality of logarithmic type and the limiting regularity condition to the harmonic heat flow.
	\newblock \emph{SIAM J. Math. Anal.} {\bf 34} (2003), no. 6, 1318--1330.  \MR{2000973}


\bibitem{RRS16}
J. C. Robinson, J. L. Rodrigo, W. Sadowski.
	\newblock The  three-dimensional Navier-Stokes equations. Classical theory. 
	\newblock  Cambridge Studies in Advanced Mathematics, 157. \emph{Cambridge University Press, Cambridge}, 2016. xiv+471 pp.  \MR{3616490}



\bibitem{Tan23}
W. Tan.
	\newblock Existence and regularity of solutions for a 3D coupled   parabolic-elliptic equations related to magnetic relaxation.
	\newblock \emph{J. Math. Anal. Appl.} {\bf 519} (2023), no. 1, Paper No. 126735, 31 pp. \MR{4493179}










\end{thebibliography}

\end{document}